%% file: article.tex
\documentclass[12pt,letterpaper]{amsart} 

\input{RVSVpreamble.tex}
\begin{document}
\pagestyle{plain}
\title{\Large{Tschirnhausen bundles of covers of the projective line}}
\author{Ravi Vakil}
\address{Dept. of Mathematics, Stanford University, Stanford CA~94305--2125}
\email{rvakil@stanford.edu}
\author{Sameera Vemulapalli}
\address{Dept. of Mathematics, Harvard University, Cambridge MA~02138--2901}
\email{vemulapalli@math.harvard.edu}
\date{September 5, 2025.}
\begin{abstract}
A degree $d$ genus $g$ cover of the complex projective line by a smooth curve $C$ yields a vector bundle  on the projective line by pushforward of the structure sheaf.  Which bundles are possible? Equivalently, which $\PP^{d-2}$-bundles over $\PP^1$ contain such covers?  (In the language of many previous papers:  what are the scrollar invariants of the cover?)

We give a complete answer in degree $4$, which exhibits the expected pathologies. We describe a polytope (one per degree)  which we propose gives the complete answer for primitive covers, i.e. covers that don't factor through a subcover. We show that all such bundles (for primitive covers) lie in this polytope, and that a ``positive proportion'' of the polytope arises from smooth covers. Moreover, we show the necessity of the primitivity assumption. Finally, we show that the  image of the map from the Hurwitz space of smooth covers to the space of bundles is not preserved by generization (for $d>4$ and $g \gg_d 1$).  
 \end{abstract}
\maketitle
\begin{figure}[ht]
\centering
\includegraphics[scale=0.25]{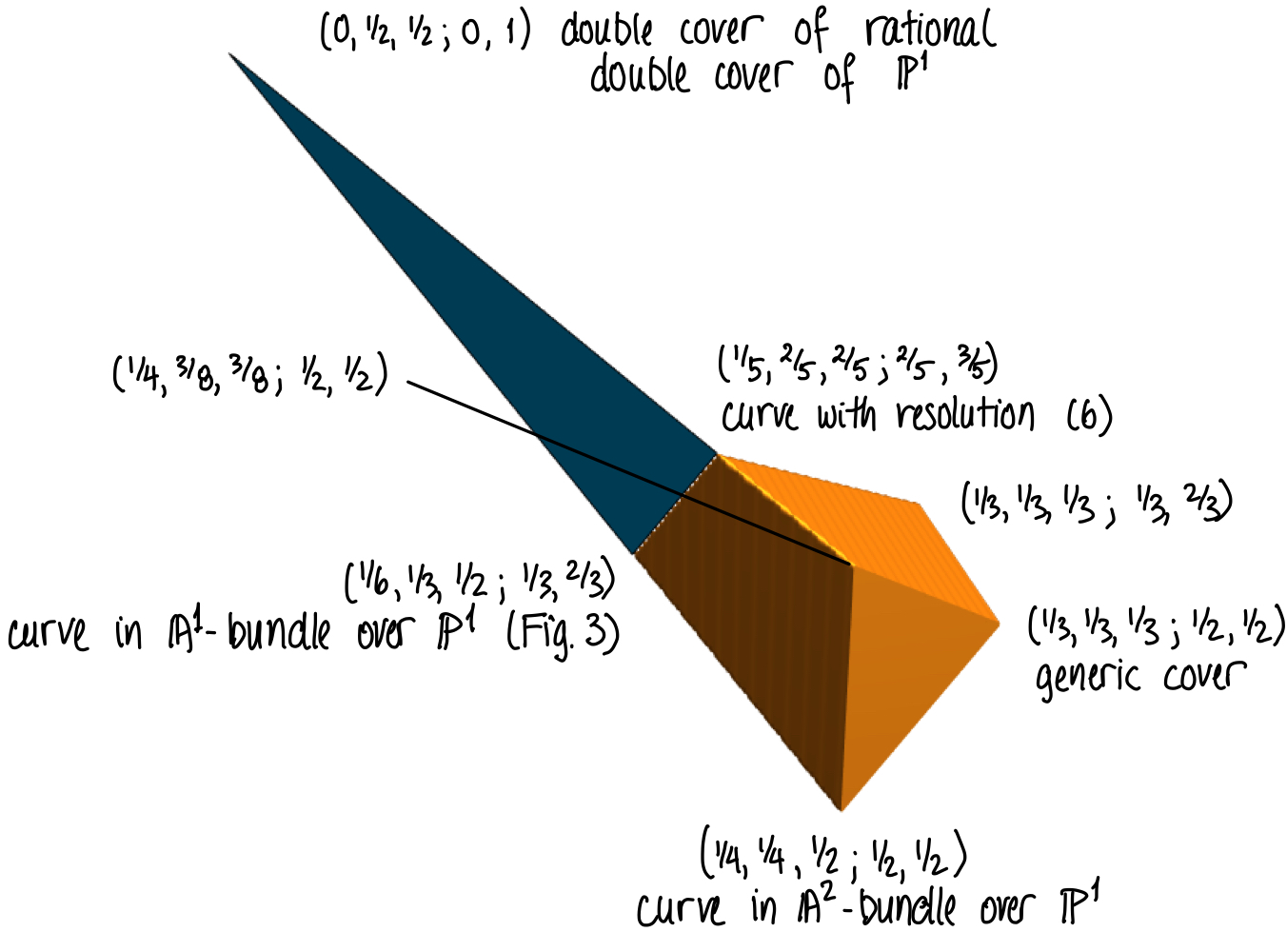}
\caption{The geography of degree $4$ covers of $\proj^1$.  More precisely  (see \S \ref{r:bigpicture}), the locus $(\ove_1, \ove_2, \ove_3; \; \ovf_1, \ovf_2)$ corresponding to tetragonal covers and their cubic resolvents:  a three-dimensional locus (mostly primitive) and a two-dimensional locus (all imprimitive).  Figure~\ref{f:triangle} is a projection of Figure~\ref{f:weirdo}.  Figure~\ref{f:monogenic} corresponds to a point of Figure~\ref{f:weirdo}.  Figure~\ref{f:3dprint} is a subset of Figure~\ref{f:weirdo}.}\label{f:weirdo}\end{figure}

\tableofcontents

{\parskip=12pt 

\section{Introduction}

Suppose  $\pi: C \rightarrow \PP^1$ is a genus $g$ degree $d$ cover of the complex projective line by a smooth irreducible complex curve ($d \geq 3$).
Then $C$ canonically embeds in a $\PP^{d-2}$-bundle over $\PP^1$. Conversely, any such genus $g$ degree $d$ cover (fiberwise scheme-theoretically nondegenerately) embedded in a $\PP^{d-2}$-bundle must be canonically embedded. Moreover, the $\PP^{d-2}$-bundle in which $C$ embeds is precisely the projectivization of the Tschirnhausen bundle of $\pi$, which we introduce below. This paper deals with the question:  {\em which bundles arise as Tschirnhausen bundles}? By the discussion above, we see that this is equivalent to the question: {\em which $\PP^{d-2}$-bundles over $\PP^1$ contain smooth genus $g$ degree $d$ covers}?

\point \label{s:early}
Let $\cE_\pi \coloneq (\pi_* \oh_C / \oh_{\PP^1})^\vee$ (the {\bf Tschirnhausen bundle of $\pi$}).
The trace map $\pi_* \oh_C \rightarrow \oh_{\proj^1}$, when multiplied by $\frac{1}{d}$, splits the exact sequence
$$\xymatrix{0 \ar[r] &  \oh_{\proj^1}  \ar[r] &  \pi_* \oh_C \ar[r] &  \cE_\pi^\vee \ar[r] &  0,}$$
so  we have $\pi_* \oh_C \overset \sim \longleftrightarrow \cE_\pi^\vee \oplus \oh_{\PP^1}$.
Vector bundles on $\proj^1$ split completely (uniquely up to non-unique isomorphism), so $\cE_\pi \cong \oh(e_1)\oplus \cdots \oplus \oh(e_{d-1})$ with $e_i$ nondecreasing.
We call the integers  $(e_1, \dots, e_{d-1})$ the {\bf scrollar invariants of the cover}.   (Caution:  in a few sources, the definition of ``scrollar invariants'' differs from these by $2$.)  
From $h^0(C, \oh_C) = h^0(\proj^1, \pi_* \oh_C)=1$, we have that $h^0(\proj^1, \cE_\pi^\vee)=0$, so the $e_i$ are all positive.  From $h^1(C, \oh_C) = g$, we have $\sum_{i=1}^{d-1} e_i = d+g-1$. 
Then the canonical embedding in question is $C \hookrightarrow \PP \cE_\pi \coloneq \operatorname{Proj}_{\PP^1} \Sym^\bullet \cE$.  (See \cite[\S 3]{lvw} for details.)  

\bpoint{Central Question}  {\em 
What can the bundle $\pi_* \oh_C$ (or equivalently, the Tschirnhausen bundle) be? Equivalently:   What are the possible $(e_1, \dots, e_{d-1})$?  If $\Hur_{d,g}$ is the stack of smooth degree $d$ genus $g$ covers,  $\Bun_{\proj^1}$ is the stack  of vector bundles on $\proj^1$, and the ``Tschirnhausen morphism''
\begin{equation}
\label{eq:defbe}\Tsch: \Hur_{d,g} \rightarrow \Bun_{\proj^1}\end{equation}
is defined by $[\pi] \mapsto  \cE_\pi$, what is the image of $\Tsch$?} \label{cq}

The answer to this question when $d \leq 4$ is given in \S \ref{s:3}--\ref{s:quartic}, and the case $d=4$ exhibits the  subtleties that should be expected in general.

For general $d$, the second author proposes an answer to  Central Question~\ref{cq} when  $\pi$ is primitive. 
(In keeping with the language for number fields,  we say that  $\pi$ is a {\bf primitive} cover if $\pi$ does not factor nontrivially through an intermediate smooth cover $\pi: C \rightarrow C' \rightarrow \proj^1$.   For example, if $d$ is prime, or if the monodromy group of $\pi$ is the full symmetric group $S_d$, then $\pi$ is necessarily primitive.) 
Define the polytope 
\begin{equation}\label{eq:V}
\polytope_d := \{ (\overline{e}_1,\dots,\overline{e}_{d-1}) \;  : \; \sum \overline{e}_i = 1,  \;  0 \leq \overline{e}_1 \leq \cdots \leq \overline{e}_{d-1},
\text{ and }\overline{e}_{i + j} \leq \overline{e}_i + \overline{e}_j \} \subset \R^{d-1}.
\end{equation}

\tpoint{Conjecture (Vemulapalli)}  {\em 
For each $d$ and $g$,  the  possible scrollar invariants for primitive degree $d$ genus $g$ smooth covers are precisely those that,  scaled down by $d+g-1$, lie in $\polytope_d$.} \label{c:sameera}

The cases $d \leq 4$ are verified in \S \ref{s:3}--\ref{s:quartic}.  
In \S \ref{s:inV}, we show that  the scrollar invariants of every primitive smooth cover 
must lie in $\polytope_d$: 

\tpoint{Theorem} {\em Suppose $\pi: C \rightarrow \proj^1$ is a primitive degree $d$ morphism from a genus $g$ curve.  Then the tuple of scrollar invariants of $\pi$, scaled by  $d+g-1$, lies in the polytope $\polytope_d$. } \label{t:one}

We first remark that the set of limit points of the scrollar invariants of primitive (not necessarily smooth) covers, scaled down to lie in $\polytope_d$, is all of $\polytope_d$; this follows from an appropriate combination of \cite[Thm.~1.1.7, Thm~1.1.9]{sameera-thesis}. Thus the  difficulty of proving Conjecture~\ref{c:sameera} lies in showing the existence of \emph{smooth} curves with prescribed scrollar invariants.

Now observe that if a smooth cover $\pi: C \rightarrow \proj^1$ has scrollar invariants not lying in $\polytope_d$, then it {\em must} factor through {\em some} intermediate cover! This is one additional way in which the scrollar invariants contain a good deal of information about the cover $\pi: C \rightarrow \proj^1$. 

We further conjecture that the answer to Question~\ref{cq} in general (not requiring primitivity) is that for each $d$, there is a finite union of rational polytopes $\polytope'_d \subset \R^{d-1}$, all of dimension $d - 2$, with the property that the possible scrollar invariants are precisely those which lie in $\polytope'_d$ when scaled down by $d+g-1$ and $\ove_1 > 0$. The example of $d=4$ (Figure~\ref{f:triangle}) shows that   $\polytope'_d$ need not equal $\polytope_d$.  Thus  the answer to Central Question~\ref{cq} for all smooth covers may, and in fact is expected to, be different from the answer to Central Question~\ref{cq} for primitive covers if $d$ is not prime! (From \cite[Thm.~1.46]{sameera-successive-minima}, one can extract an expectation of $\polytope'_d$ if $d$ is $12$, or the product of two primes, or a prime power.  A precise conjecture for general $d$ seems intractable, however.)

Moreover, if $d$ is prime, Conjecture~\ref{c:sameera} implies that the image of $\Tsch$ is convex (interpreted as a subset of $\Z^{d-2}$).  In \S \ref{s:six}, we show that the image of $\Tsch$ is not convex for $d=6$ and sufficiently large genus (so $d=6$ is the first ``nonconvex degree'', if Conjecture~\ref{c:sameera} holds for $d=5$), and we indicate how to show the analogous statement for any $d$ that is not a prime power.  

In \S \ref{s:positive}, we  show that a ``positive proportion'' of the polytope  $\polytope_d$ is actually achieved.  More precisely, 
we show the following.
 We say a sequence of real numbers $(x_0,\dots,x_d)$ is {\bf concave} if $x_{i + 1} - x_i \leq x_{i} - x_{i - 1}$ for $1 \leq i \leq d-1$.  

\tpoint{Theorem}
\label{mainconstructionthm}
{\em Choose integers $(e_1,\dots,e_{d-1}) \in \ZZ_{\geq 1}^{d-1}$ such that $(0,e_1,\dots,e_{d-1},0)$ is a concave sequence. 
Then there is a primitive cover $\pi: C \rightarrow \proj^1$ by a smooth (irreducible) curve with $\cE_\pi \cong \oplus_{i=1}^{d-1} \oh(e_i)$.}

We produce curves with ``concave'' scrollar invariants by constructing certain singular curves on Hirzebruch surfaces. The crucial input is Wood's explicit and striking parametrization of curves on Hirzebruch surfaces using binary forms. (Wood's work builds on a long history of associating rings to binary forms, including \cite{na}, \cite{df}, \cite{bm}, and \cite{dcds}.) Wood's construction gives us enough control over the algebra to control the scrollar invariants of the normalization. 

\epoint{Remarks} (a)
The ``concave locus'' in $\polytope_d$ (where the sequence $(\ove_1, \dots, \ove_{d-1})$ is nondecreasing and $(0,\ove_1, \dots, \ove_{d-1},0)$ is a concave sequence) is a full-dimensional subpolytope of $\polytope_d$.  (For example, in $\polytope_4$ as shown in Figure~\ref{f:triangle}, it is the ``upper half triangle'' that is the convex hull of $(1/6, 1/3, 1/2)$, $(1/5, 2/5, 2/5)$, and $(1/3, 1/3, 1/3)$.)  

(b)  But note further that the $e_i$ in Theorem~\ref{mainconstructionthm} are not required to be nondecreasing, so by considering the $e_i$ produced by this theorem, and reordering them to be nondecreasing, we get a union of  full-dimensional polytopes in $\polytope_d$ that are shown to be realizable by our construction.   (For example, in $\polytope_4$, by taking $e_2 \geq e_3 \geq e_1$ in Theorem~\ref{mainconstructionthm} and then reordering, we obtain the ``right half triangle'' that is the convex hull of $(1/5, 2/5, 2/5)$, $(1/3, 1/3, 1/3)$ and $(1/4, 1/4, 1/2)$; see Figure~\ref{f:triangle}.  So together with the previous remark, this theorem shows that three quarters of $\polytope_4$ is ``realizable''.)  

(c) The number of such polytopes is  $2^{d-2}$:  Given an ordered list, there are $2^{d-1}$ ways of rearranging the list elements so they are monotonically increasing then monotonically decreasing, and Theorem~\ref{mainconstructionthm} shows the existence of covers in each of these cases.  Reversing the sequence gives the same construction.  (The case $d=4$ will make this combinatorial argument clear to the reader.)     

(d)  A variation of Theorem~\ref{mainconstructionthm} can  be used to get the final quarter of $\polytope_4$ (omitted because all of $\polytope_4$ is already shown to be achieved in \S \ref{s:quartic}), but the authors do not see how to get all of $\polytope_d$ in general. See Remark~\ref{r:variant} for more details.

Question~\ref{cq} has been considered in a number of geometric contexts.  We mention here other work showing that scrollar invariants in $\polytope_d$ are achieved. The main theorem   of \cite{coppens}  shows that a particular type of scrollar invariant is achieved, and is the first result of this sort of which we are aware for general degree $d$.    \cite[Ex.\ 1.3.7]{ll} also describes an interesting case. \cite[\S 7]{cvz} shows that certain highly imbalanced scrollar invariants are achieved, albeit not for smooth curves.  
Most dramatically, the main theorems of \cite{ballico}, \cite{coppens2} (a simpler second proof of the result of \cite{ballico}), and \cite{dp} can be interpreted as showing that in some asymptotic sense, all points ``near'' the most ``generic'' or ``balanced'' point of $\polytope_d$ (where the $\ove_i$ are equal) are  achieved, and these include cases not covered by what we show in \S \ref{s:positive}.  (We did not check to see if \cite{dp} gives a region of $\polytope_d$ of positive volume, because the statement of their result does not immediately translate in this way.  The results of \cite{ballico} and \cite{coppens2} do not give regions of positive volume.)  See also the discussion of \cite{ko} in \S \ref{s:quartic2}.

In \S \ref{s:notflat} we ask if the image of the morphism $\Tsch: \Hur_{d,g} \rightarrow \Bun_{\proj^1}$ (of \eqref{eq:defbe}) is preserved by generization, i.e. if $\cE$ is the Tschirnhausen bundle of a smooth cover, and $\cE'$ is a generization of $\cE$, is it true that $\cE'$ is \emph{also} the Tschirnhausen bundle of a smooth cover? There are a number of reasons to hope the image of $\Tsch$ is preserved by generization. If so, then the simple deformation theory of $\Bun_{\proj^1}$ would allow us to prove the existence of curves with other scrollar invariants.  It would also allow us to prove statements about fibers of $\Tsch$ by looking at the most degenerate cases, which are far more tractable (the ``minimal scrollar invariants'' correspond to curves on a Hirzebruch surface not meeting the directrix). The fact that the image of $\Tsch$ seems to be of such a nice form (conjecturally a polytope in the primitive case) might further give hope of good behavior.  

And indeed the image of $\Tsch$ is preserved by generization  when $2g + 1 \leq d$ or $d \leq 4$. Surprisingly, we show (in Proposition~\ref{p:notclosed}) that the image of $\Tsch: \Hur_{d,g} \rightarrow \Bun_{\proj^1}$ is not preserved by generization ($d \geq 5$, $g \gg_d 1$). In the case where $d$ is not a prime power, one might expect this from the nonconvexity of the image of $\Tsch$, but when $d$ is a prime power, we expect the image to be convex and yet it is not preserved by generization. In the primitive case, we even show that the image of the restriction $\Tsch^p: \Hur_{d,g}^p \rightarrow \Bun_{\proj^1}$ to primitive covers is not preserved by generization ($d \geq 4$, $g \gg_d 1$).

Question~\ref{cq} is related to questions about the lattice structure of rings of integers in number fields (see \S \ref{s:arithmetic}); scrollar invariants are naturally analogous, when appropriately scaled, to successive minima of rings of integers. However, we may ask a more refined question: we are interested in not just \emph{which} Tschirnhausen bundles occur, but also how often they occur, and this question has a number theoretic analogue as well. In \S \ref{s:arithmetic} we discuss how the latter question can be described via a conjectural density function on $\polytope_d$ in both the arithmetic and geometric settings. We then further speculate that these density functions extend to the Grothendieck ring.

\bpoint{Central remaining question}  It remains to either show that the remaining points of $\polytope_d$ are achieved by primitive smooth covers, or (more dramatically) to show that there is some further obstruction to certain Tschirnhausen bundles being achieved for primitive smooth covers. 

\bpoint{Conventions}    We work over $\C$ for simplicity, although the results readily extend in obvious ways to other fields. Henceforth, the phrase {\bf  cover of $\proj^1$} means  a finite morphism $C \rightarrow \proj^1$ from an irreducible  curve $C$, and
 the phrase {\bf smooth cover of $\proj^1$} means further that $C$ is smooth. 
  Throughout $d$ will refer to the degree of $\pi$, and $g$ will refer to the arithmetic genus of $C$.  Throughout this article, we fix a splitting 
\begin{equation}\label{eq:fixsplitting}
\pi_* \oh_C \overset \sim \longleftrightarrow \oh_{\proj^1} \oplus \oh_{\proj^1}(-e_1) \oplus \cdots \oplus \oh_{\proj^1}(-e_{d-1})
\end{equation}
where the $e_i$ are as described.  It will be convenient to take the convention $e_0=0$, so \begin{equation}
\label{eq:fixsplitting2}
\pi_* \oh_C = \oplus_{i=0}^{d-1} \oh_{\proj^1}(-e_i).
\end{equation}

We follow common notational conventions for Hirzebruch surfaces $\F_a \rightarrow \proj^1$.  In $\Pic \F_a$, let $F$ be the fiber class, let $E$ be the directrix (which is the unique section of negative self-intersection if 
$a>0$, and has self-intersection $-a$). A section of $\FF_a$  not intersecting the directrix has class $D=E+ a F$.  (See  for example \cite[\S V.2]{hartshorne} for standard facts about Hirzebruch surfaces.)

\bpoint{Acknowledgments}  We thank Manjul Bhargava, Hendrik Lenstra, and Hannah Larson for many useful conversations. Melanie Wood's work has also been an inspiration. The second author was supported by the NSF under grant number DMS2303211. We would also like to warmly thank Edoardo Ballico, Marc Coppens, Anand Deopurkar, Aaron Landesman, Hannah Larson, Hendrik Lenstra, Daniel Litt, Anand Patel, and Will Sawin for helpful comments on an earlier version of this draft. Finally, the referee's careful reading significantly improved the manuscript.

\section{Warm-up:  Low degree ($d \leq 3$) and low genus ($2g+1 \leq d$)}

\label{s:3}

\bpoint{The trivial cases $d \leq 2$}  The rank and degree of the bundle $\cE_\pi$ determine the bundle for low-rank bundles.  We quickly observe that for $d=1$ we have $\cE_\pi = 0$, and for $d=2$ we have $\cE_\pi = \oh_{\proj^1}(g+1)$.  

\bpoint{The straightforward case $d=3$}
\label{ss:cubic}
Before getting to the truly enlightening case of $d=4$, we consider the $d=3$ case.    Consider a degree $3$ cover $\pi: C \rightarrow \proj^1$ from a genus $g$ curve.  Then $\pi_* \oh_C = \oh_{\proj^1 } \oplus \oh_{\proj^1 } (-e_1)\oplus \oh_{\proj^1 } (-e_2)$
where $e_2 \geq e_1 \geq 1$, and $e_1+e_2 = g+2$.   Let $m=e_2-e_1 $ for convenience. The canonical embedding 
$C \hookrightarrow \PP \cE$ embeds $C$ in the Hirzebruch surface 
$\FF_m \rightarrow \PP^1$.  A smooth irreducible trisection of $\FF_m$ is in the class $3D +kF$ for some $k \geq 0$, and a short calculation shows that such a curve has genus $3m -2+ 2k$.  Now for $C$ to be such a curve, we must have $3m-2+2k = g$, from which
$3(e_2-e_1)-2+2k = e_1+e_2-2$, from which $k=2e_1-e_2$.
As $k \geq 0$, we have $2e_1 \geq e_2$.  We conclude that for 
trigonal curves, we have $e_1 +e_2 = g+2$, and $2 e_1 \geq e_2 \geq e_1$.  (This is clearly $\polytope_3$ as defined in \eqref{eq:V}.)   By choosing general curves in the appropriate linear system 
$\oh_{\FF_m}(3D+kF)$, we see that there is indeed a trigonal cover 
$\pi: C \rightarrow \proj^1$ for each of these values of $(e_1, e_2)$.  We have thus completely answered the Central Question in this case.

We see that the scrollar invariants, suitably normalized, are precisely those that lie on the line segment  joining $(1/2, 1/2)$ and $(1/3, 2/3)$. 

\epoint{Aside} The natural follow-up question is: as $g \rightarrow \infty$, how \emph{often} do various points on this line segment occur from smooth curves. In \S \ref{s:density}, we define a density function $\rho^{geo}$ from the aforementioned line segment to the real numbers, which quantifies the ``likelihood'' that a curve lands ``near'' a given point on the line segment. It turns out that the curves are not ``equidistributed'' on this line segment as $g \rightarrow \infty$ (equivalently, $\rho^{geo}$ is not a constant function, but indeed is a nontrivial linear function on the line segment). 

We further conjecture that this fits into a more general framework (see \S \ref{s:density} for more details). We define an arithmetic analogue $\rho^{arith}$ capturing information about the shapes of rings of integers in cubic fields, as a lattice. Moreover, we define (conjecturally) a measure $\pi^{gr}_{ref}$ on the line segment taking values in the Grothendieck ring of stacks, appropriately localized. We then define a (conjectural) function $\rho^{gr}_{ref}$, also taking values in the (appropriately completed) Grothendieck ring of stacks, which captures the point masses of this measure. We conjecture that $\dim \circ \rho^{gr}_{ref} = \rho^{geo} = \rho^{arith}$.

\bpoint{Low genus}    If the genus is low compared to the degree, the Tschirnhausen bundle is forced to be balanced (``as close to equal as possible given that they have to be integers''), as is made precise by the following. 

\tpoint{Proposition} \label{p:triangle} {\em If $g \leq \frac {d-1} 2$, then 
$\pi_* \oh_C = \oh_{\proj^1} \oplus \oh_{\proj^1}(-1)^{\oplus (d-1-g)} \oplus \oh_{\proj^1}(-2)^{\oplus g}.$}

The proof showcases an idea, illustrated in the following lemma, that will be used in the proof of Theorem~\ref{t:one}. There are a number of essentially analogous lemmas in the arithmetic case; see \cite[Lem.~6.4]{pikertrosen} or the proof of \cite[Thm.~1.6]{2torsion}.

\tpoint{Lemma} \label{l:square}  {\em Suppose $C$ is an irreducible curve, and $\pi: C \rightarrow \proj^1$ is a degree $d$ cover.  Then (following the notation of \eqref{eq:fixsplitting2}) $e_i + e_{d-1-i} \geq e_{d-1}$ ($0 \leq i \leq d-1$).}

\bpf  
Choose a point $\infty \in \proj^1$, and choose a coordinate $t$ on $\AAA^1 = \proj^1 \setminus \{ \infty \}$, i.e., an isomorphism $\proj^1 \setminus \{ \infty \}  \cong \Spec \CC[t]$.  Then the splitting $\pi_* \oh_C = \oh_{\proj^1} \oplus \oh_{\proj^1}(-e_1) \oplus \cdots \oplus \oh_{\proj^1}(-e_{d-1})$
induces a splitting of the $\CC[t]$-algebra
$\Ga(\AAA^1, \pi_*(\oh_C))$ into
$\Ga(\AAA^1, \pi_*(\oh_C)) = \C[t] \oplus \CC[t] x_1 \oplus \cdots \oplus \CC[t] x_{d-1}$
(as a $\CC[t]$-module); here we have chosen a generator  $x_i$ of the $i$th summand  of 
$\Ga(\AAA^1, \pi_*(\oh_C))$.
Now $1$, $x_1$, \dots, $x_{d-1}$ are a basis for $K(C)$ as a $K(\proj^1)$-module (where $K(\cdot)$ indicates the function field) and the restriction of the multiplication in the $\oh_{\proj^1}$-algebra 
$$
\pi_* \oh_C \otimes \pi_* \oh_C \rightarrow \pi_* \oh_C$$ is simply multiplication $K(C) \otimes_{K(\proj^1)} K(C) \rightarrow K(C)$ in the field. Let $M$ be the $d \times d$ matrix given by taking the field multiplication map and projecting to the last coordinate. Observe that $M$ is invertible because it represents a perfect pairing. Indeed, given any nonzero $v_1 \in K(C)$, there is $v_2 \in K(C)$ such that $v_1v_2 = x_{d-1}$; simply take $v_2 = x_{d-1}/v_1$.

Now notice that for every $0 \leq a,b < d$, the summand corresponding to \begin{equation}\oh(-e_a) \otimes \oh(-e_b) \rightarrow \oh(-e_{d-1}) \label{eq:alg}\end{equation}
is an element of $H^0( \oh_{\proj^1}(e_a+e_b-e_{d-1}))$, and is thus $0$ if $e_a+e_{b}<e_{d-1}$. But if $e_i+e_{d-1-i}<e_{d-1}$, then for $a \leq i$, $b \leq d-i$, the summand \eqref{eq:alg} is zero, so $M$ cannot be invertible, as it has an $(i+1) \times (d-i)$ submatrix that is zero.
\epf

\epoint{Proof of Proposition~\ref{p:triangle}} 
Recall that  $1 \leq e_1 \leq \cdots \leq e_{d-1}$ and $\sum e_i = d-1+g$.  If $g=0$, the result is then immediate, so assume $g>0$. 
If $(e_1, \dots, e_{d-1})$ is  balanced, then $e_1=\cdots = e_{d-g-1}=1$, and $e_{d-g} = \cdots = e_{d-1}=2$.  
A quick thought will show that if $(e_1, \dots, e_{d-1})$ is {\em not} balanced, then $e_{d-1}>2$ and $e_{d-g}=1$.  But then $e_{g-1} \leq e_{d-g}=1$ gives $e_{g-1}+e_{d-g} < e_{d-1}$, contradicting Lemma~\ref{l:square}. \epf

\section{Key Example:  Quartic covers ($d=4$)}

\label{s:quartic} We now turn to the case $d=4$.

\tpoint{Theorem} \label{t:quartic} {\em  The possible scrollar invariants of tetragonal covers are precisely those that, when suitably scaled, are contained in the closed triangle shown in Figure~\ref{f:triangle}, minus the point $(0, 1/2, 1/2)$. Every such triple is realized by a smooth bi-hyperelliptic cover. Moreover, the possible scrollar invariants of primitive tetragonal covers are those that, when suitably scaled, lie in the closed polytope $\polytope_4$.} \label{t:tet}

\begin{figure}[ht]
\includegraphics[scale=0.20]{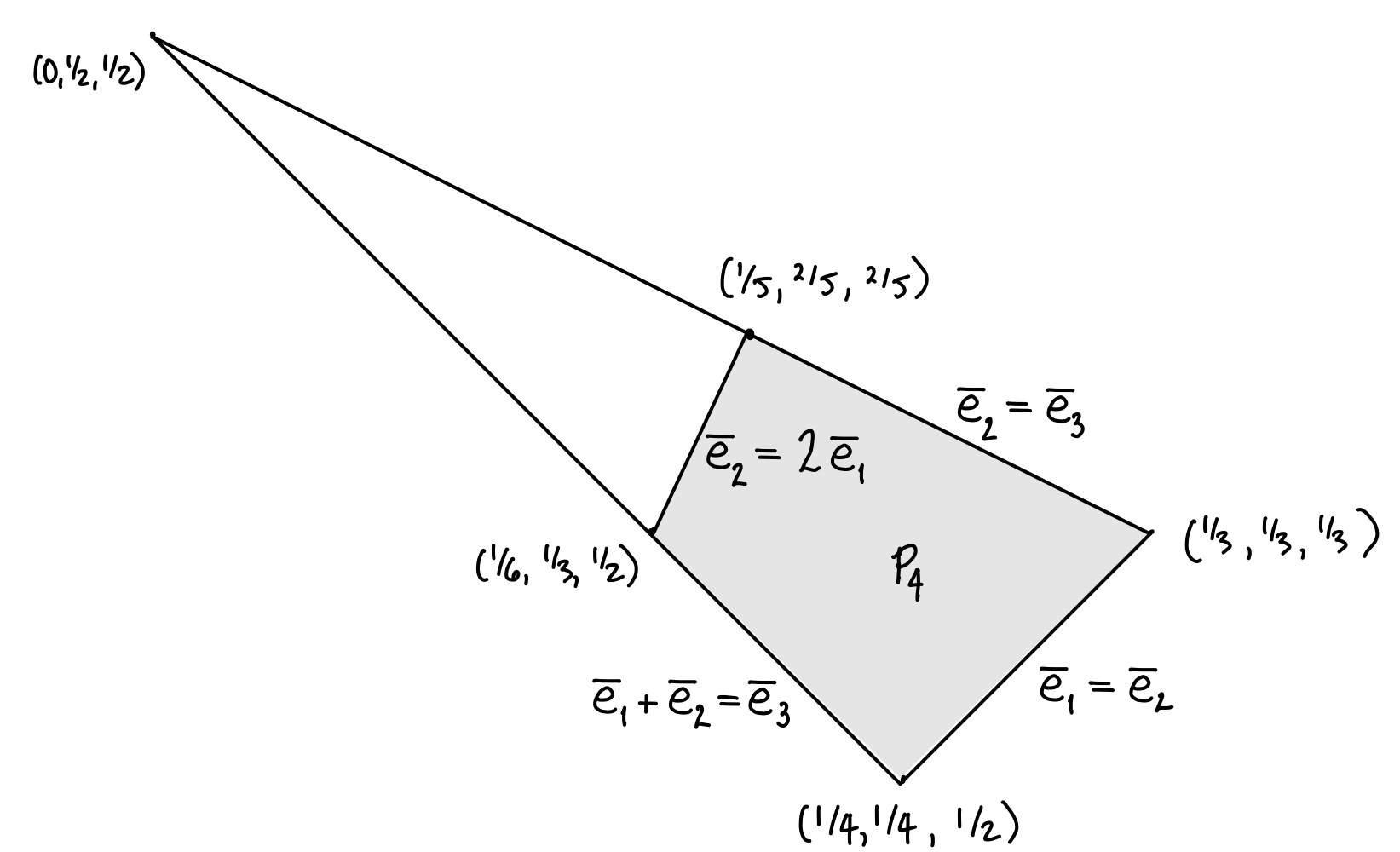}
\caption{Scrollar invariants of tetragonal covers (showing points $(\overline{e}_1, \overline{e}_2, \overline{e}_3) \in \R^3$):  the big triangle is $\polytope'_4$,  the shaded quadrilateral is $\polytope_4$, and the small triangle minus its intersection with $\polytope_4$ necessarily corresponds to imprimitive covers.  }\label{f:triangle}
\end{figure}

\epoint{History} \label{s:quartic2}The case of quartic
covers has a long and rich history.  We cannot hope to describe the
highlights in a short paragraph, but we mention a few relevant results.
We point out in particular \cite{CaL, pa, sameera-successive-minima} for constraints on scrollar invariants of tetragonal covers. Also see \cite{cm2} for linear series on tetragonal covers (and \cite{cm1} for arbitrary degree);
and \cite{casnatidc} and \cite{casnati} on degree $4$ covers of more general varieties.  Finally, strikingly:  \cite{ko} proves that every scrollar invariant
that (suitably scaled, and with $e_1 > 0$) lies in $\polytope'_4$ is achieved by a tetragonal cover, but the covers described there are all {\em imprimitive}.

\epoint{Observations} Before proving Theorem~\ref{t:tet}, we make some observations.  

Any tetragonal cover $\pi: C \rightarrow \proj^1$ with $\pi_*\oh_C = \oh \oplus \oh(-e_1) \oplus \oh(-e_2) \oplus \oh(-e_3)$ with $e_2 > 2e_1$ necessarily is not primitive; we will see that it factors through a hyperelliptic cover (and such curves have Galois/monodromy group smaller than $S_4$, and the generic such curve has monodromy group $D_4$).

Some of the vertices of $\polytope_4$ correspond to curves with (other) recognizable geometric meaning.  For example, the generic  genus $g$ tetragonal curve has balanced scrollar invariants, so when $3 \mid g$, we have $\cE \cong \oh ( \frac 1 3 (g+3) )^{\oplus 3}$.
Thus the vertex $(1/3, 1/3, 1/3)$ corresponds to ``generic covers''.

Second, curves $C$ of class $4D$ on the Hirzebruch surface $\FF_n$ have genus $g=6n-3$.  By embedding $\FF_n$ (interpreted as a 
$\proj^1$-bundle over $\proj^1$) as a conic bundle in the corresponding $\proj^2$-bundle over $\proj^1$, we find curves 
embedded in $\proj_{\proj^1} ( \oh \oplus \oh(n) \oplus \oh(2n))$, see Figure~\ref{f:monogenic}.  Suitably normalizing, this 
corresponds to $\overline{e}_1$, $\overline{e}_2$, and $\overline{e}_3$ in arithmetic progression; this gives the vertex $(1/6, 1/3, 1/2)$.  (These are the ``least 
generic'' scrollar invariants for primitive covers.) Conversely, curves lying at the vertex $(1/6,1/3,1/2)$ are precisely those of class $4D$ on a Hirzebruch surface. These might be called monogenic covers, as they lie in $\AAA^1$-bundles over $\proj^1$. Another reason for the term ``monogenic'' is that in this case, $\oh \oplus \oh(-e_1)$ generates $\pi_* \oh_C$ as an $\oh_{\proj^1}$-algebra.

\begin{figure}[ht]
\includegraphics[scale=0.14]{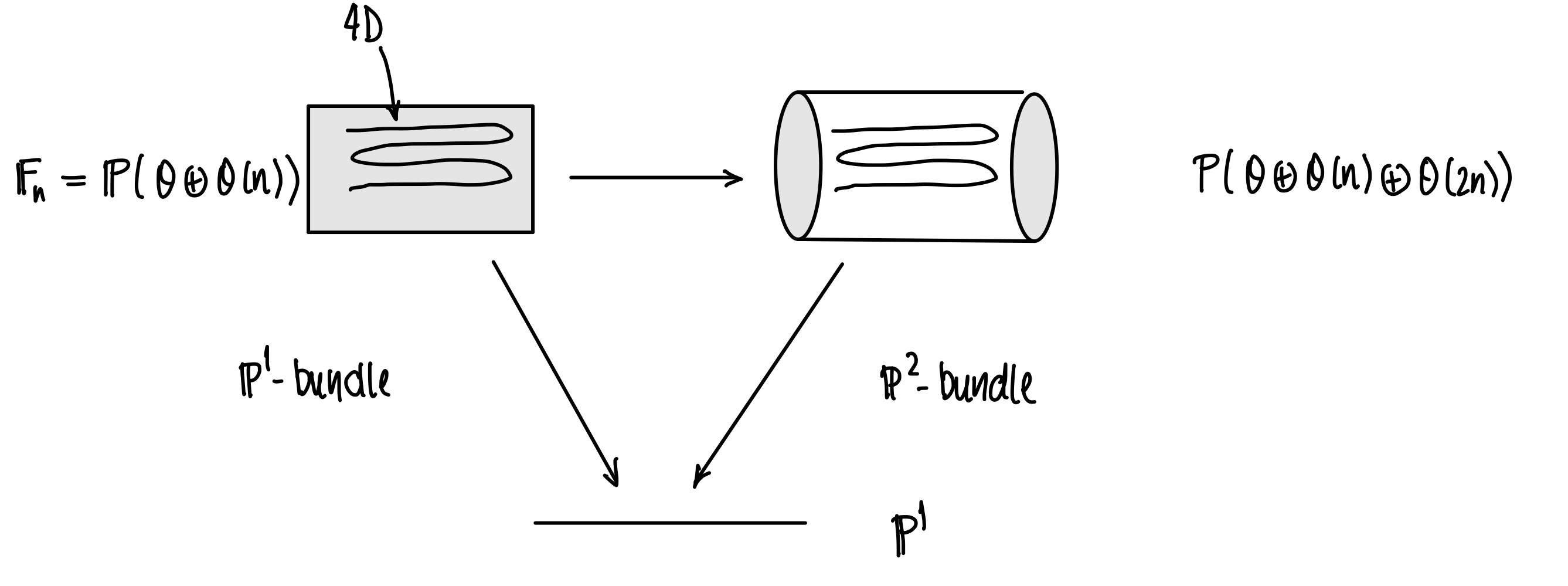}
\caption{Monogenic covers:  $(\ove_1, \ove_2, \ove_3) = (1/6, 1/3, 1/2)$}\label{f:monogenic}
\end{figure}

The vertex $(1/4, 1/4, 1/2)$ corresponds similarly to curves in $\AAA^2$-bundles over $\proj^1$. The curves lying on this vertex are precisely those that do not intersect the Hirzebruch surface $\Proj(\oh(-e_1) \oplus \oh(-e_2)) = v(x_3)$ in their embedding into $\PP(\cE_{\pi}^{\vee})$.

The vertex  $(1/5, 2/5, 2/5)$ is harder to describe. These curves have the property that there exists a Hirzebruch surface $\FF_{e_3 - e_1}$ given by the vanishing of the linear equation $ax_2 + bx_3$, where $a,b \in \CC$ and $x_1,x_2,x_3$ are the coordinates of $\PP(\cE_{\pi}^{\vee})$. The defining linear equation is obtained as follows: there is a distinguished conic in the conic bundle, and it contains the line $v(x_2,x_3)$. The tangent space to that line forms the Hirzebruch surface $\FF_{e_3 - e_1}$.    Algebraically, if $\cE^\vee = \oh_{\proj^1}(e_1) \oplus \oh_{\proj^1}(e_2) \oplus \oh_{\proj^1}(e_3) = \oh_{\proj^1}(k) \oplus \oh_{\proj^1}(2k) \oplus \oh_{\proj^1}(2k)$, then $\xymatrix{C  \ar@{^(->}[r] &  \proj \cE^\vee  \ar[r]^\pi & \proj^1}$ with resolution
\begin{equation}\label{eq:five}
\xymatrix{0 \ar[r] & \pi^* \oh_{\proj^1}(-5k)  \otimes \oh_{\proj \cE^\vee}(-4) \ar[r] &
\pi^* \left( \oh_{\proj^1}(-2k) \oplus \oh_{\proj^1}(-3k)    \right) \otimes \oh_{\proj \cE^\vee}(-2)   \ar[r]  & \\
& \oh_{\proj \cE^\vee} \ar[r] 
& \oh_C \ar[r] & 0.
}\end{equation}
The numerology  (which applies similarly to the other vertices) is explained by \eqref{eq:3dprint}.

\epoint{Proof of Theorem~\ref{t:tet}}

We begin by reviewing more geometry about tetragonal covers.  
(This is classical, so rather than trying to find the first reference, we  refer the reader to \cite[\S 3]{lvw}.) Suppose $\pi: C \rightarrow \proj^1$ is such a cover.    Then $C \hookrightarrow \proj \cE_\pi^\vee$ is the complete intersection of two divisors $\Delta_1$ and $\Delta_2$  in $\proj \cE_\pi^\vee$. The divisors $\Delta_1$ and $\Delta_2$ are of
class $2H + f_1$ and $2H + f_2$, where $H$ is the class of $\oh_{\PP E^{\vee}}(1)$, and \begin{equation}
\label{eq:linear}f_1 + f_2 = d+g-1 = e_1 +e_2 + e_3.\end{equation}  Say $f_2 \geq f_1$ without loss of generality.  The ``equation'' of $\Delta_k$ is of the following form.  Let $x_1$, $x_2$, $x_3$ be the projective coordinates on the $\proj^2$-bundle $\proj \cE_\pi^\vee$  corresponding to $e_1$, $e_2$, $e_3$ (so $x_i$ has ``weight'' $-e_i$).   Let $(s,t)$ be the projective coordinates on the base $\proj^1$.  Then   $\Delta_k$ has equation
$$\xymatrix@=0pt{
a^k_{11}(s,t) x_1^2 & + & a^k_{12}(s,t) x_1 x_2 & + & a^k_{22}(s,t) x_2^2 \\
+ & a^k_{13}(s,t) x_1 x_3 & + & a^k_{23}(s,t) x_2 x_3 & &  = 0. \\
& + & a^k_{33}(s,t) x_3^2
}$$
Here $a^k_{ij}(s,t)$ is homogeneous (in $s$ and $t$) of degree $e_i + e_j - f_k$.  In particular, the equation for $\Delta_k$ has ``weight'' $-f_k$.  Notice that  if $e_i + e_j -f_k<0$, then necessarily $a^k_{ij}(s,t)=0$.

Rather than looking for which triples $(e_1,e_2,e_3)$ are achievable from such $\pi$, we ask which quintuples $(e_1, e_2, e_3; f_1, f_2)$ are achievable.  

By our conventions, we have $e_1 \leq e_2 \leq e_3$ and $f_1 \leq f_2$. 

\point \label{s11}
Now if $f_1> 2 e_1$, then $f_2>2 e_1$ as well, so  $a^1_{11}=a^2_{11}=0$, and the equations for $\Delta_1$ and $\Delta_2$ have no $x_1^2$ term, so both $\Delta_1$ and $\Delta_2$ contain the section $V(x_2, x_3)$ of $\proj \cE_\pi^\vee$.  But then the complete intersection is reducible, not an irreducible quadruple cover.  Hence we have $f_1 \leq 2 e_1$.   

\point \label{s22}
Also, if  $f_2> 2 e_2$, then $f_2> e_1 +  e_2$ and $f_2 > 2e_1$ as well, so  $a^2_{11} = a^2_{12}=a^2_{22}=0$, and hence the equation for $\Delta_2$ is divisible by $x_3$, and $\Delta_2$ is thus reducible, so the complete intersection $\Delta_1 \cap \Delta_2$ cannot be an integral curve.  Hence we have $f_2 \leq 2 e_2$.  

\point
The five inequalities cut out a polytope $\qpolytope'$ with vertices 
\begin{eqnarray}\label{eq:qprime}
(\ove_1, \ove_2, \ove_3; \; \ovf_1, \ovf_2) &=& \left( \frac 1 3, \frac 1 3, \frac 1 3;  \; \frac 1 2, \frac 1 2 \right), \; \left( \frac 1 4, \frac 1 4, \frac 1 2; \; \frac 1 2, \frac 1 2 \right), \; \left( \frac 1 4, \frac 3 8, \frac 3 8; \; \frac 1 2, \frac 1 2 \right), \\ & & 
\left( \frac 1 3, \frac 1 3, \frac 1 3; \; \frac 1 3, \frac 2 3 \right), \;  \left( 0, \frac 1 2, \frac 1 2;  \; 0, 1 \right) \nonumber
\end{eqnarray}
where $\ovf_k \coloneq f_k/(g+3)$.  (Compare the first three coordinates of each vertex with Figure~\ref{f:triangle}.)

\point \label{s123}Notice further that if $f_2> e_1+e_3$, then $f_2> e_1+e_2$ and $f_2>2 e_1$, so $a^2_{13} = a^2_{12} = a^2_{11}=0$, so the conic bundle $\Delta_2$ is of the form 
\begin{equation}\label{eq:hyp}a_{22}^2 x_2^2 + a^2_{23} x_2 x_3 + a^2_{33} x_3^2=0,\end{equation} so the fibers (over geometric points of $\proj^1$) are pairs of lines (both passing through $x_2=x_3=0$), and a short argument shows that the complete intersection factors through a hyperelliptic cover (with equation given by \eqref{eq:hyp}, except taken in the $\proj^1$-bundle with projective coordinates $x_2$ and $x_3$).    Thus 
for primitive covers, we have an additional linear inequality $f_2 \leq e_1 + e_3$. 

\point This, along with the five previous inequalities, cut out a convex polytope $\qpolytope$ with vertices
\begin{eqnarray} \label{eq:3dprint}
(\ove_1, \ove_2, \ove_3; \; \ovf_1, \ovf_2) &=& \left( \frac 1 3, \frac 1 3, \frac 1 3;  \; \frac 1 2, \frac 1 2 \right), \; \left( \frac 1 4, \frac 1 4, \frac 1 2; \; \frac 1 2, \frac 1 2 \right), \; \left( \frac 1 4, \frac 3 8, \frac 3 8; \; \frac 1 2, \frac 1 2 \right), \\ & & 
\left( \frac 1 3, \frac 1 3, \frac 1 3; \; \frac 1 3, \frac 2 3 \right), \;  \left( \frac 1 5, \frac 2 5, \frac 2 5; \;  \frac 2 5, \frac 3 5 \right), \;  \left( \frac 1 6, \frac 1 3, \frac 1 2;  \; \frac 1 3, \frac 2 3  \right) \nonumber
\end{eqnarray}
See Figure~\ref{f:3dprint} for a picture of $\qpolytope$. 
The astute reader will observe that we have proved Theorem~\ref{t:one} in the case $d=4$:  the projection of the polytope  $\qpolytope$ to the first three coordinates is $\polytope_4$.  We have shown a necessary condition for a quintuple to come from a primitive tetragonal cover, and we now show sufficiency.

\begin{figure}[ht]
\centering
\includegraphics[scale=0.20]{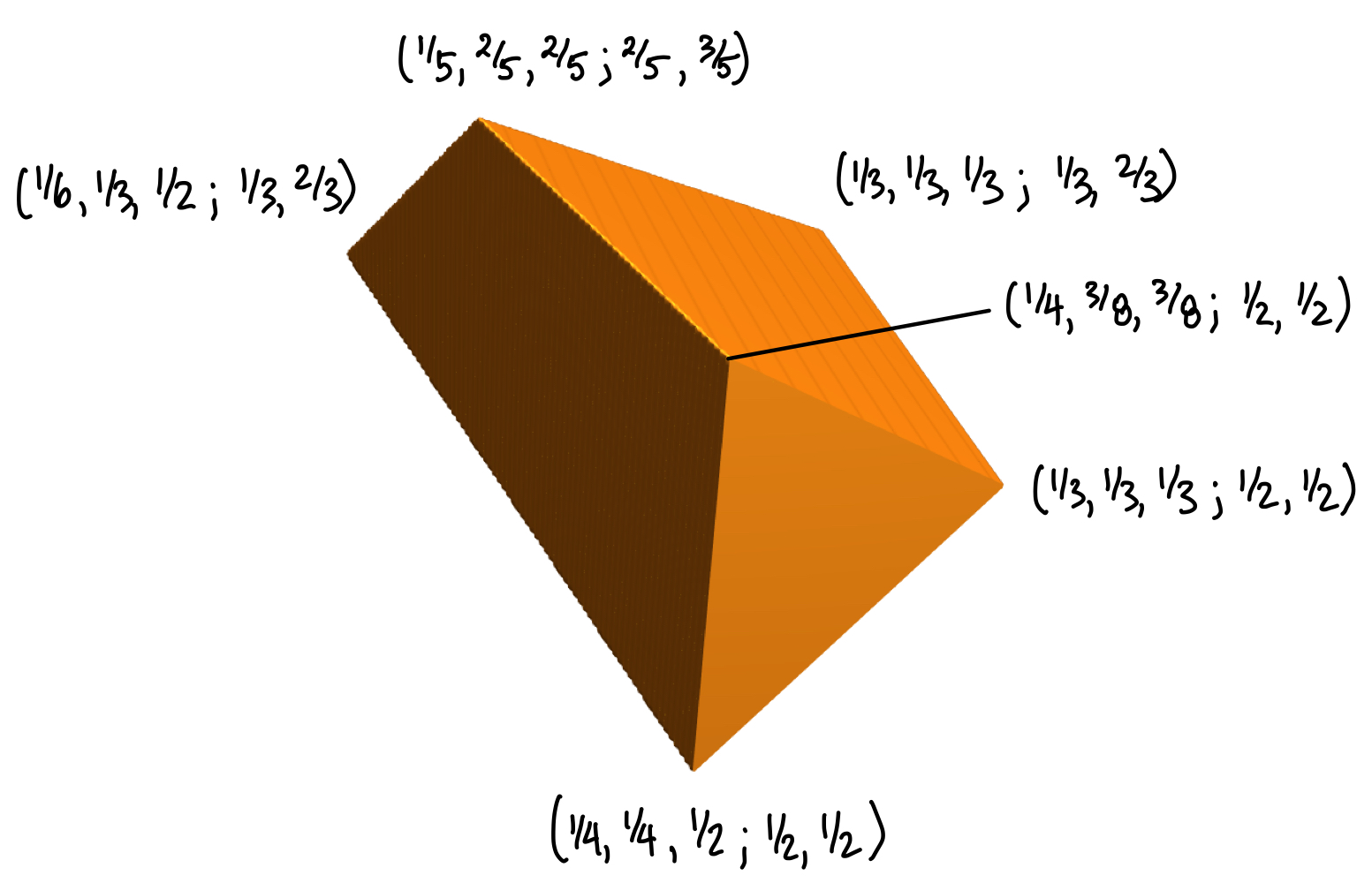}
\caption{The polytope showing possible $(\ove_1, \ove_2, \ove_3; \;  \ovf_1, \ovf_2)$  for primitive tetragonal covers, portrayable in $\R^3$ using the linear conditions $\ove_1+\ove_2+\ove_3 = \ovf_1+\ovf_2=1$ (cf.\ Figure~\ref{f:weirdo}).
}\label{f:3dprint}
\end{figure}

We observe next that if $2 e_1-f_1>0$, then the divisor $\Delta_1$ is ample: it is easily seen to separate points and tangent vectors on the $\proj^2$-bundle $\proj \cE^\vee_\pi$.   This condition holds for all of $\qpolytope$ except for a quadrilateral $Q$ on the boundary, so we assume now that $\Delta_1$ is ample.  Fix now a point $(e_1, e_2,e_3; f_1, f_2) \in \Z^5$ arising from a point in $\qpolytope \setminus Q$, by scaling the latter point by $g+3$.  

By Bertini's Theorem, if we can show that a generally chosen $\Delta_2$ is smooth away from finitely many points, then ampleness of the class of $\Delta_1$ means that a generally chosen $\Delta_1 \cap \Delta_2$ will give a smooth curve (in particular a quartic cover of $\proj^1$).  Now a general fiber (over a generally chosen point of $\proj^1$ of a generally chosen $\Delta_2$ will have the equation $? x_1 x_3 + ? x_3^2 + ? x_2 x_3 + ? x_2^2 + ??$, with the single question marks generally chosen nonzero values, and the double question marks indicating that possibly there are additional monomials $x_1^2$ and $x_1 x_2$ appearing.  Thus the general fiber is easily seen to be a smooth conic, so the generally chosen $\Delta_2$ is generically smooth over $\proj^1$, as desired.

\epoint{The remaining case $2e_1=f_1$ (when $\Delta_1$ is not ample)}
\label{s:D1notample}

The edge case where $2e_1 =f_1$ is similar.  In this case, $\Delta_1$ does not separate points:  it contracts the points on the section $x_2=x_3=0$.  But away from this section, the divisor class $\Delta_1$  separates points and tangent vectors, so $\Delta_1$ is the pullback of a very ample divisor contracting $\proj \cE^\vee_\pi$ along this section.  The previous arguments then apply essentially without change.

\bpoint{Existence of primitive covers}
\label{s:existenceofprimitive}
At this point we have shown that for all points of $\qpolytope$, the corresponding general complete intersection is a tetragonal cover of $\proj^1$ by a smooth curve.  We need to show next that this cover is primitive, i.e., that it does not factor through a double cover of $\proj^1$.    If the general such cover is imprimitive, then the monodromy group of the cover lies in $D_4$, so the resolvent cubic is reducible. To show this is not the case, we need to exhibit a single example where the resolvent cubic is irreducible.  We do this by exhibiting an example where the resolvent cubic is irreducible over the generic point of $\proj^1 = \Proj \C[s,t]$.    

In preparation, notice that for all $$(i,j,k) \in   \{ (1,1,1), (2,3,1), (1,3,2), (2,2,2) \}, $$
we have $e_i+e_j-f_k \geq 0$
($2e_1 \geq f_1$ by \S \ref{s11}, then  $e_2+e_3 \geq f_1$  from $e_3 \geq e_2 \geq e_1$, $2e_2 \geq f_2$ by \S \ref{s22}, $e_1+e_3 \geq f_2$ by \S \ref{s123}).
Furthermore, we cannot have $e_2+e_3-f_1 =0$, as then $2e_1-f_1=0$ as well, from which $e_1=e_2=e_3$, but then $2 e_1 = f_1 \leq (f_1+f_2)/2=(e_1+e_2+e_3)/2 = 3e_1/2$ from which $e_1=0$, contradicting the irreducibility of $C$.

Consider  now the  special case where $\Delta_1$ is given by
$s^{2 e_1-f_1} x_1^2 + 2   s^{e_2+e_3-f_1-1} t  x_2 x_3=0$, and $\Delta_2$ is given by $2 s^{e_1 + e_3 - f_2} x_1 x_3 +   s^{2 e_2 - f_2}  x_2^2=0$.  

The cubic resolvent is (readily shown to be, and well-known to  be)
$$
\det \left(  x \left( \begin{array}{ccc} s^{2 e_1-f_1} & 0 & 0 \\ 0 & 0  &  s^{e_2+e_3-f_1-1} t  \\ 0 &  s^{e_2+e_3-f_1-1} t & 0  \end{array}\right) + y \left(  \begin{array}{ccc} 0  & 0 &  s^{e_1 + e_3 - f_2}  \\ 0 & s^{2 e_2 - f_2} & 0  \\s^{e_1 + e_3 - f_2} & 0 & 0  \end{array}\right)   \right)
$$
which is $$-  s^{2 e_1 + 2 e_2 + 2 e_3 }  \left(  (t/s)^2 (x/s^{f_1})^3  +  (y/s^{f_2})^3 \right) $$
which is easily seen to be an  irreducible cubic  in $\C(s,t)  [x,y]$ (for example using Eisenstein's criterion). 

\bpoint{Completing the proof of Theorem~\ref{t:tet}}
The following claim will prove that the possible scrollar invariants of primitive tetragonal covers of genus $g$ are those in the (closed) polytope $\polytope_4$, scaled by $g+3$.  This will complete the proof of the part of Theorem~\ref{t:tet} dealing with primitive covers.

\tpoint{Claim} {\em Suppose $0 < e_1 \leq e_2 \leq e_3$ are integers with $\frac 1 {g+3} (e_1, e_2, e_3) \in \polytope_4$.  Then there are integers $f_1 \leq f_2$ with $\frac 1 {g+3} (e_1, e_2, e_3; f_1, f_2) \in \qpolytope$.}
\bpf
Take $f_1 = \min(2 e_1, \lfloor \frac {e_1 + e_2 + e_3} 2 \rfloor )$, and $f_2 = g+3-f_1$.  We show that $(e_1, e_2, e_3; f_1, f_2) \in \qpolytope$ by showing that it satisfies the six desired inequalities.  Four of them are automatic, so we need only check (a)  $f_2 \leq 2 e_2$ and (b) $f_2 \leq e_1 + e_3$.

For (a):  $f_2 \leq 2 e_2$ is equivalent to $e_1 + e_2 + e_3 - f_1 \leq  2 e _2$, which is equivalent to $e_1 + e_3 \leq e_2 + f_1$. So it suffices to show that (a1)  $e_1 + e_3 \leq e_2 + 2 e_1$ and  (a2) $e_1 + e_3 \leq e_2 +\lfloor \frac {e_1 + e_2 + e_3} 2 \rfloor$.   (a1)  is one of the inequalities defining $\polytope_4$.   (a2)  is  easily manipulated to  $e_2 \geq \lceil \frac {e_1 + e_2 + e_3} 4 \rceil$.  But by inspection of the vertices of $\polytope_4$, $\ove_2 \geq 1/4$, so $e_2 \geq \frac { e_1 + e_2 + e_3} 4$, and $e_2 \in \Z$, so this case is complete. 

Inequality (b) is equivalent to $e_2 \leq f_1$, so it suffices to show (b1) $e_2 \leq 2 e_1$ and (b2) $e_2 \leq \lfloor \frac {e_1 + e_2 + e_3} 2 \rfloor$.  (b1) is one of the inequalities defining $\polytope_4$.  
By inspection of the vertices of $\polytope_4$, $\ove_2 \leq 1/2$, so $e_2 \leq \frac {e_1 + e_2 + e_3} 2$, so (b2) follows as $e_2 \in \Z$. \epf 

Finally, note that \cite{ko} proves that if $(\ove_1,\ove_2,\ove_3)$ lies in $\polytope_4'$, then it is realized by some smooth bi-hyperelliptic cover. In all cases, the smoothness of $C$ implies the irreducibility of $C$, as $h^0(C,\oh_C) =  \oplus_{i = 0}^3 h^0(\PP^1,\oh(-e_i)) = 1$. This concludes the proof of Theorem~\ref{t:tet}. \epf

\bpoint{The remainder of $\qpolytope'$: the exclusively ``bi-hyperelliptic'' or ``$D_4$ monodromy'' locus} 
\label{fillingq'}We now consider the remaining points in the polytope described by \eqref{eq:qprime}. We are in the regime where $f_2> e_1 + e_3$ so $\Delta_2$ is given by \eqref{eq:hyp}.    This 
surface is singular along the entire section $x_2=x_3=0$.  Thus any complete intersection $\Delta_1 \cap \Delta_2$ will be singular {\em unless} $\Delta_1$ does not intersect this section $x_2=x_3=0$.  Translation:  unless $a^1_{11}(s,t)$ is a constant, i.e., $f_1 = 2 e_1$, there cannot be a complete intersection that is a smooth curve.  Thus outside of the locus $\qpolytope$ all remaining smooth curves are in the plane $f_1 = 2 e_1$.     In this case, $\Delta_1$ is {\em not} ample, but its failure of ampleness is our friend:  a general section of $\oh(\Delta_1)$ will be nonzero along all of $x_2=x_3=0$, and away from this locus Bertini's Theorem applies.    Then in this case the previous arguments apply without change, and a general such complete intersection will give a tetragonal cover by a regular curve, except of course this cover factors through a hyperelliptic cover of $\proj^1$.
Note only that the point $(0,1/2, 1/2)$ is excluded because from the start we had a strict inequality $e_1>0$. 

\bpoint{Quintuples in $\qpolytope$ arising from smooth bi-hyperelliptic covers} We now ask which quintuples $(e_1,e_2,e_3;f_1,f_2)$ occur from smooth bi-hyperelliptic covers. If the quintuple $(\ove_1,\ove_2,\ove_3;\ovf_1,\ovf_2)$ is contained in $\qpolytope' \setminus \qpolytope$, then the discussion above provides an affirmative answer precisely when $f_1 = 2e_1$. We now prove that any quintuple $(\ove_1,\ove_2,\ove_3;\ovf_1,\ovf_2) \in \qpolytope$ occurs from a smooth bi-hyperelliptic curve.

If $2e_1 = f_1$, then we have shown in \S \ref{s:existenceofprimitive} that the general corresponding complete intersection is primitive (hence not bi-hyperelliptic). Therefore, to show the existence of a smooth bi-hyperelliptic curve with the prescribed invariants, we restrict to a special subfamily of complete intersections. Namely, take the conic bundle $\Delta_2$ to be general subject to the condition that $a_{11}^2 = a_{12}^2 = a_{13}^2 = 0$. Now, the argument of  \S \ref{fillingq'} applies without change, and a general such complete intersection is a smooth bi-hyperelliptic curve (a short calculation shows that the cubic resolvent is reducible, and hence the curve is bi-hyperelliptic).

Now, if $2e_1 > f_1$, again we have shown in \S \ref{s:existenceofprimitive} that the general corresponding complete intersection is primitive, so we restrict to a special family. Choose a general polynomial $p(s,t)$ of degree $e_3 - e_2$. Consider conic bundles $\Delta_1$ and $\Delta_2$ subject to the conditions that $a_{11}^2 = a_{12}^2 = a_{12}^1 = 0$ and $a_{23}^2 = pa_{22}^2$ and $a_{23}^1 = pa_{22}^1$. Fix a general choice of $\Delta_2$. Observe that it is smooth away from finitely many points; a general fiber is a smooth conic and no fiber over $\proj^1$ is a double line.  Again, the cubic resolvent is easily seen to be reducible. 

Now, notice that the class of $\Delta_1$ separates points away from $V(a_{22}^2)$. It separates tangent vectors away from the union of $V(x_2,x_3)$ and $V(x_1,x_3)$. Therefore by Bertini's theorem, a general intersection $\Delta_1 \cap \Delta_2$ is smooth on $\PP(\cE)^{\vee} \setminus \{V(a_{22}^2) \cup V(x_1,x_3) \cup V(x_2,x_3)\}$. 

To complete the proof, observe that for a general choice of $\Delta_1$ and $\Delta_2$, we have:
\begin{enumerate}
    \item there are no singular points of $\Delta_1 \cap \Delta_2$ contained on  $V(a_{22}^2)$ (on these fibers, the intersection of $\Delta_1 \cap \Delta_2$ is generically $4$ distinct points);
    \item the intersection of $\Delta_1 \cap \Delta_2 \cap V(x_1,x_3)$ is trivial (this follows from the fact that $a_{22}^2$ and $a_{22}^1$ share no common factors);
    \item and there are no singular points of $\Delta_1 \cap \Delta_2$ contained on  $V(x_2,x_3)$ (to see this, compute the rank of the Jacobian matrix at any point of intersection).
\end{enumerate}

\bpoint{Quintuples $(e_1,e_2,e_3;f_1,f_2)$ arising from smooth covers}
 \label{r:bigpicture}For a number of reasons it is worth considering which quintuples $(e_1, e_2, e_3; f_1, f_2)$ correspond to smooth  covers.  (For example:  these track the entire resolution of $C \hookrightarrow \proj \cE^\vee$; the $f_i$ correspond to the resolvent cubic; and more.)  The above discussion is enough to determine the achievable quintuples, and the result is more striking than when we just identify the possible $(e_1, e_2, e_3)$; we get the union of two polytopes of different dimensions (see Figure~\ref{f:weirdo}).

\bpoint{Fun remarks about the ``fin''}
There are even more phenomena worth pointing out.  For the purposes of this discussion, we call the two-dimensional triangular offshoot the ``fin'' (where to be precise, the fin excludes the boundary where it meets the $3$-dimensional polytope shown in Figure~\ref{f:3dprint}). The moduli space of quartic covers of fixed genus is irreducible -- so perhaps it is surprising to see that its image (admittedly, under $\Tsch$, and asymptotically as $g \rightarrow \infty$) has two pieces. 

\epoint{The ``bulk'' of quartic covers live on the ``balanced'' vertex of $\qpolytope$}
The generic quartic cover of fixed genus is balanced (i.e. $e_3 -e_1 \leq 1$). This implies that the ``bulk'' of the ``mass'' of the moduli space is concentrated on the ``balanced'' vertex, where the $\ove_i$ are all equal and the $\ovf_j$ are both equal. 

In \S \ref{s:density} we make this statement precise by defining a \emph{density function} $\rho_{ref}^{\geo} \colon \R^5 \rightarrow \R$ which sends a point $(\ove_1,\ove_2,\ove_3;\ovf_1,\ovf_2)$ to a quantity which reflects the ``proportion'' of quartic curves with those given invariants. Surprisingly, $\rho_{ref}^{\geo}$ is a piecewise linear function, supported on the region pictured on Figure~\ref{f:weirdo}, and continuous on its support. We give an explicit formula for $\rho_{ref}^{\geo}$ (Proposition~\ref{p:quartic-rho}), and show that $\rho_{ref}^{geo}(1/3,1/3,1/3;1/2,1/2) = 1$, and this is where it obtains its maximum value. 

\epoint{But, a codimension $2$ sublocus of $\cH_{4,g}$ lives as far away from the balanced vertex as possible}

\begin{figure}[ht]
\centering
\includegraphics[scale=0.35]{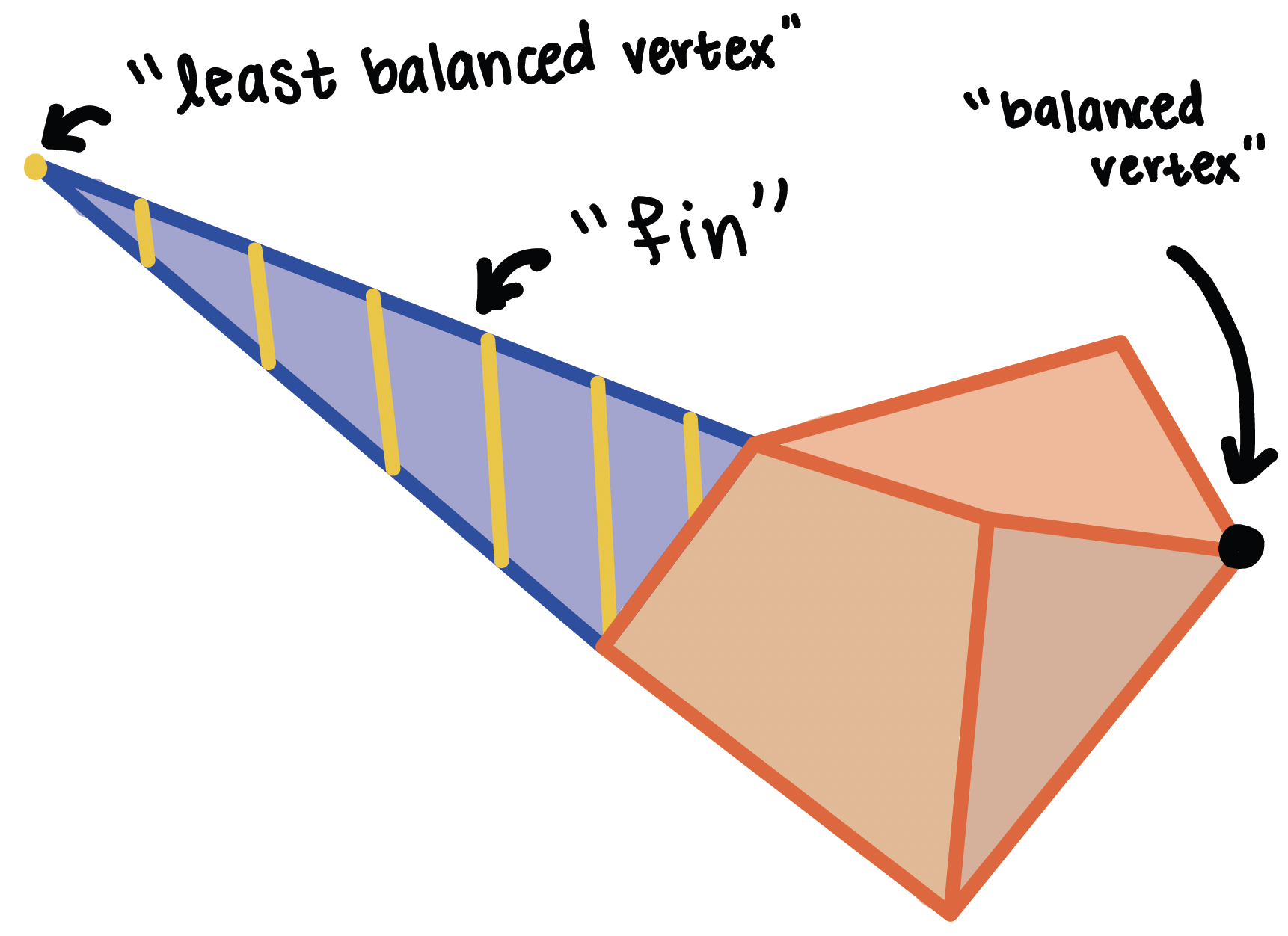}
\caption{A depiction of the geography of degree $4$ covers of $\proj^1$. The ``fin'' is drawn in blue, and the yellow lines indicate the loci where the hyperelliptic subcover has fixed genus. The leftmost yellow point corresponds to the hyperelliptic subcover having genus $0$.}\label{f:hlines}\end{figure}

However, there is one other point where $\rho^{geo}_{ref}$ obtains its maximum value of $1$, and that is at the {\em least} balanced vertex, $(0, 1/2, 1/2; 0,1)$, which is the point farthest from $\polytope_4$. Here, the picture is as follows. We have already shown that covers in the fin are bi-hyperelliptic, but in fact we can do better: such covers have a \emph{distinguished} quadratic subcover, and the genus of that quadratic subcover is precisely $e_1 - 1$.

\bpoint{Lemma} { \em
Let $C \rightarrow \PP^1$ be a quartic cover whose scrollar invariants lie on the fin. Then $C$ has a distinguished quadratic subcover, and that quadratic subcover has genus $h = e_1 - 1$.
}

\bpf 
Because $C$ lies in the fin, $e_2 > e_1$. Consider the multiplication map $\pi_* \oh_C \otimes \pi_*\oh_C \rightarrow \pi_* \oh_C$ and fix a splitting of $\pi_* \oh_C$. Let $\cA$ be the subbundle $\oh_{\proj^1} \oplus \oh_{\proj^1}(-e_1)$. Because $e_2 > 2e_1$, the image of $\cA \otimes \cA$ under the multiplication map lies in $\cA$; in other words, $\cA$ is a quadratic $\oh_{\proj^1}$-subalgebra of $\pi_* \oh_C$. Therefore $C \rightarrow \proj^1$ factors through the  hyperelliptic curve $\underline{\Spec} \;  \cA$, which visibly has genus $e_1 - 1$. 
\epf

Upon fixing a genus $g$, each possible value of $h$ gives rise to a line in the fin (explicitly, the line $\frac{h}{g+3} = \ove_1 - \frac{1}{g+3}$). Upon fixing $h$, the ``bulk'' of the moduli space of quartic covers factoring through genus $h$ hyperelliptic curves lies where the line meets the edge of the fin given by $\ove_2= \ove_3$ (a glance at Figure~\ref{f:hlines} may be enlightening). Concretely: fix a value $0 \leq \overline{h} \leq 1/4$. Then, the maximum value of $\rho_{ref}^{geo}$ restricted to the hyperplane $\ove_1 = \overline{h}$ is achieved precisely when $\ove_2= \ove_3$. 

Now, the maximum value of $\rho_{ref}^{geo}$ on the fin is attained at the least balanced vertex $(0,1/2,1/2;0,1)$, which is precisely the point where $\overline{h} = 0$. In this case, the dimension of the moduli space of bi-hyperelliptic quartic covers whose quadratic subcover is a genus $0$ curve is precisely $2g + 1$. Compare this to the dimension of the moduli space of quartic covers, which is $2g + 3$. The two spaces have comparable dimension as $g \rightarrow \infty$, and the difference in dimensions is always $2$! This explains why $\rho_{ref}^{geo}$ attains the same value at both vertices. 

\epoint{Arithmetic meaning of the fin}
This discussion is related to the remarkable arithmetic fact due to Bhargava that, when ordered by absolute discriminant, a positive portion of quartic number fields have Galois group $D_4$, not $S_4$; see \cite[Thm.~4]{bh}. 

Qualitatively, our discussion above shows that there are lots of quartic covers whose scrollar invariants lie either at the ``balanced'' vertex or the ``least balanced'' vertex. The scrollar invariants of curves are analogous to the successive minima of rings of integers in number fields (see \S \ref{s:density} for details). It is possible to show \cite[Thm.~5.5]{sameera3} that when ordered by absolute discriminant, that $S_4$-fields tend to lie at the ``balanced'' vertex, and $D_4$-fields tend to lie at the ``least balanced'' vertex. The fact that the density function $\rho_{ref}^{geo}$ takes value $1$ on both vertices reflects the fact that a positive portion of quartic number fields have Galois group $D_4$.

\begin{figure}[ht]
\centering
\includegraphics[scale=0.35]{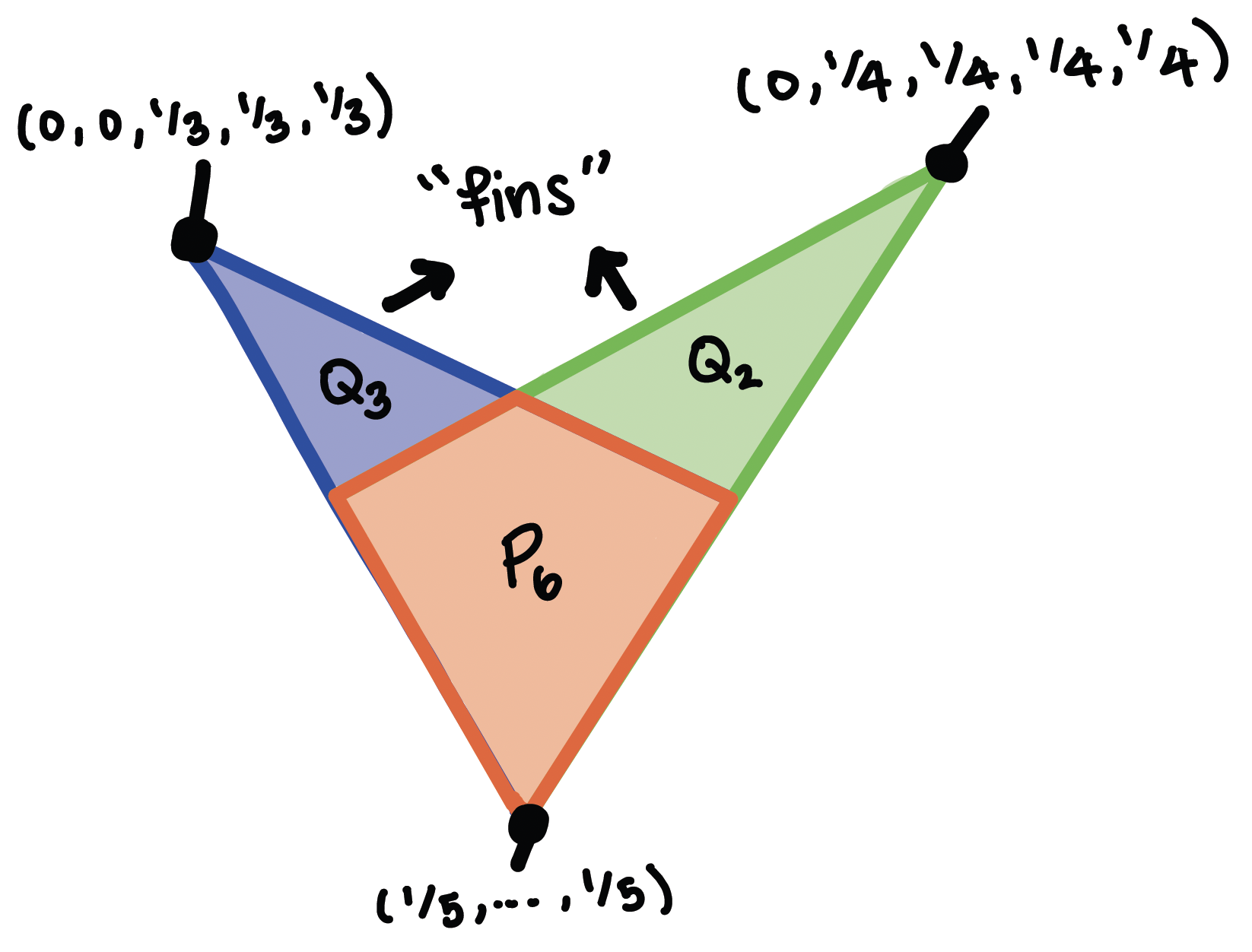}
\caption{A tale of two fins: in this cartoon drawing, $Q_3$ is the union of the blue and orange region, $Q_2$ is the union of the green and orange region, and $\polytope_6$ is their intersection. The two unbalanced points at the top, as well as the balanced point at the bottom, are labeled.}\label{f:twofins}\end{figure}

\epoint{``Fins'' occur in all composite degrees}
This is not a peculiarity of $d=4$; there is a similar story with ``fins'' in all non-prime degrees.  (This should be seen as analogous to the fact that Malle's conjecture guesses that there are many field extensions with Galois group not $S_d$ when $d$ is composite). 

The key ideas are visible in the case $d = 6$. Define

\[
Q_3 := 
\left.
\begin{cases}
(\overline{e}_1,\dots,\overline{e}_{5}) \;  : \; \sum \overline{e}_i = 1,  \;  0 \leq \overline{e}_1 \leq \cdots \leq \overline{e}_{5},
\text{ and } \\
\ove_2 \leq 2\ove_1, \; \ove_4 \leq \ove_3 + \ove_1, \; \ove_5 \leq \ove_3 + \ove_2, \; \text{ and } \ove_5 \leq \ove_1 + \ove_4
\end{cases} \right\}
\subset \R^{5}.
\]

\[
Q_2 := 
\left.
\begin{cases}
(\overline{e}_1,\dots,\overline{e}_{5}) \;  : \; \sum \overline{e}_i = 1,  \;  0 \leq \overline{e}_1 \leq \cdots \leq \overline{e}_{5},
\text{ and } \\
\ove_3 \leq \ove_1 + \ove_2, \; \ove_4 \leq 2\ove_2, \; \ove_5 \leq \ove_3 + \ove_2, \; \text{ and }\ove_5 \leq \ove_1 + \ove_4
\end{cases} \right\}
\subset \R^{5}.
\]

A suitable generalization of Theorem~\ref{t:one} shows that the scrollar invariants of a sextic curve lie in $Q_3 \cup Q_2$; moreover it is not hard to check that $Q_3 \cap Q_2$ is precisely $\polytope_6$. (For a proof, see \cite[Thm.~1.26]{sameera-successive-minima}). One might conjecture that the possible scrollar invariants of sextic curves are precisely those that (appropriately scaled) lie in $Q_2 \cup Q_3$, but this is not necessary for the discussion here. 

\tpoint{Proposition} { \em
If $C$ is a sextic cover whose scrollar invariants lie in $Q_2 \smallsetminus \polytope_6$, then $C$ factors through a distinguished double cover with scrollar invariant $e_1$. Similarly, if $C$ is a sextic cover whose scrollar invariants lie in $Q_3 \smallsetminus \polytope_6$, then $C$ factors through a distinguished triple cover whose scrollar invariants are $(e_1,e_2)$.
}

\bpf
Consider the multiplication map $\pi_* \oh_C \otimes \pi_*\oh_C \rightarrow \pi_* \oh_C$ and fix a splitting of $\pi_* \oh_C$. Let $\cA$ be the subbundle $\oh_{\proj^1} \oplus \oh_{\proj^1}(-e_1)$. If  $\cA$ is a quadratic $\oh_{\proj^1}$-subalgebra of $\pi_* \oh_C$, then $C \rightarrow \proj^1$ factors through the  hyperelliptic curve $\underline{\Spec} \;  \cA$, which visibly has scrollar invariant $e_1$.

If $C$ has scrollar invariants lying in $Q_2 \smallsetminus \polytope_6$, then  $e_2 > 2e_1$ or $e_4 > e_1 + e_3$. If $e_2 > 2e_1$, then the image of $\cA \otimes \cA$ under the multiplication map lies in $\cA$ so $\cA$ is a quadratic $\oh_{\proj^1}$-subalgebra of $\pi_* \oh_C$. If $e_4 > e_1 + e_3$, let $\mathcal{B}$ be the subbundle $\oh_{\proj^1} \oplus \dots \oplus \oh_{\proj^1}(-e_3)$. Then the image of $\cA \otimes \mathcal{B}$ lies in $\mathcal{B}$, so $\cA$ is a quadratic $\oh_{\proj^1}$-subalgebra of $\pi_* \oh_C$, as required.

If $C$ has scrollar invariants lying in $Q_3 \smallsetminus \polytope_6$, then  $e_3 > e_1 + e_2$ or $e_4 > 2e_2$. In the first case, letting $\cA$ be the subbundle as above and letting $\mathcal{C} = \oh_{\proj^1} \oplus \oh_{\proj^1}(-e_1) \oplus \oh_{\proj^1}(-e_2)$, we have that the image of $\cA \otimes \mathcal{C}$ lies in $\mathcal{B}$, and therefore $\mathcal{C}$ is a cubic $\oh_{\PP^1}$-algebra. If $e_4 > 2e_2$, then letting $\mathcal{B}$ be as above, we see that the image of $\mathcal{C} \otimes \mathcal{C}$ lies in $\mathcal{B}$, so $\mathcal{C}$ is a cubic $\oh_{\PP^1}$-algebra as required. (See Lemma~\ref{l:hlx} for more details.) 
\epf

Let $\rho^{geo}$ be the (conjectural) density function constructed in \S \ref{s:density}. It is conjecturally (Conjecture~\ref{conj:rho}) a piecewise linear function supported on a polytope, and we expect that its support will be $Q_3 \cup Q_2$. Because a general sextic cover has balanced scrollar invariants (i.e. $e_5 - e_1 \leq 1$), we expect that $\rho^{geo}(1/5,\dots,1/5) = 1$. (And this is indeed the analogue of the observation that the ``bulk'' of quartic covers lie on the ``balanced'' vertex). 

Upon fixing a genus $g$, each possible value of $e_1$ gives rise to a $3$-dimensional polytope in $Q_2 \smallsetminus \polytope_6$.  We now consider the dimension of the moduli space of sextic covers factoring through double covers with scrollar invariant $e_1$. A short calculation shows that this dimension is largest precisely when $e_1 = 1$, and indeed the dimension of that moduli space is $2g + 3$. A general such curve has scrollar invariants $(1,e_2,\dots,e_5)$ where $e_5 -e_2 \leq 1$. 

Similarly, each possible value of $e_1,e_2$ gives rise to a $2$-dimensional polytope in $Q_3 \smallsetminus \polytope_6$. The dimension of the moduli space of sextic covers factoring through triple covers with scrollar invariants $e_1, e_2$ is largest when $e_1 = e_2 = 1$, and in this case the dimension is also $2g + 3$. A general such curve has scrollar invariants $(1,1,e_3,e_4,e_5)$ where $e_5 - e_3 \leq 1$. 

The dimension of the moduli space of sextic covers is $2g+7$. Because the difference in dimensions $(2g+ 7) - (2g+3) = 4$ is a constant, we have that if the density function $\rho^{geo}$ is defined on the points $(0,1/4,1/4,1/4,1/4)$ and $(0,0,1/3,1/3,1/3)$, then the density function achieves its maximum value of $1$ on these points.  

In general composite degree $d$, the ``fin'' points will be of the form $(0,\dots,0,1/(d-k),\dots,1/(d-k))$ for each nontrivial proper divisor $k \mid d$. They are attained by covers of the form $C \rightarrow \PP^1 \rightarrow \PP^1$ where the first map has degree $d/k$ and the map $\PP^1 \rightarrow \PP^1$ has degree $k$.

\section{Scrollar invariants of primitive covers  lie in the polytope $\polytope_d$}

\label{s:inV}

We now prove Theorem~\ref{t:one}.
Our argument should be interpreted as  an adaptation of \cite[Thm.~1.1.8]{sameera-thesis} (see also \cite[\S 2.2]{sameera-thesis}) to a geometric context.

\epoint{Observation}
\label{o:obs}Consider the $\oh_{\proj^1}$-algebra structure on the vector bundle $\pi_* \oh_C = \oh_{\proj^1} \oplus \oh_{\proj^1}(-e_1) \oplus \cdots \oplus \oh_{\proj^1}(-e_{d-1})$ \eqref{eq:fixsplitting}. By convention we take $e_0=0$, and call the ``$\oh(-e_i)$'' summand the $i$th summand.    Under the product structure in this sheaf of algebras, the product of the $i$th summand and the $j$th summand, decomposed again into summands, must be zero in any summand $\oh(-e_k)$ where $e_k > e_i + e_j$, as $\Hom(\oh(-e_i) \otimes \oh(-e_j), \oh(-e_k)) = 0$ in that case.  

\point
\label{point:coords}
Choose a point $\infty \in \proj^1$, and choose a coordinate $t$ on $\AAA^1 = \proj^1 \setminus \{ \infty \}$, i.e., an isomorphism $\proj^1 \setminus \{ \infty \}  \cong \Spec \CC[t]$.  Then the splitting $\pi_* \oh_C = \oh_{\proj^1} \oplus \oh_{\proj^1}(-e_1) \oplus \cdots \oplus \oh_{\proj^1}(-e_{d-1})$
induces a splitting of the $\C[t]$-algebra
$\Ga(\AAA^1, \pi_*(\oh_C))$ into
$\Ga(\AAA^1, \pi_*(\oh_C)) = \C[t] \oplus \CC[t] x_1 \oplus \cdots \oplus \CC[t] x_{d-1}$
(as a $\CC[t]$-module); here we have chosen a generator  $x_i$ of the $i$th summand  of 
$\Ga(\AAA^1, \pi_*(\oh_C))$.
Considered as a rational section of $\pi_* \oh_C$, $x_i$ has a pole at $\infty$ of order $e_i$.  
Now $1$, $x_1$, \dots, $x_{d-1}$ are a basis for $K(C)$ as a $K(\proj^1)$-module (where $K(\cdot)$ indicates the function field).  

\point 
We now recall a beautiful result of 
Hou, Leung and Xiang \cite{hlx}, which is a field-theoretic analogue of Kneser's theorem in additive combinatorics.  (\cite{bsz} generalized the result, but we do not need that version here.)

Suppose $L/K$ is a finite degree field extension of characteristic zero fields.   (We won't need this fact, but we remark that the results stated here in fact hold over any field, \cite{bsz}.)  Suppose that  $A$ and $B$ are two nonzero $K$-subvector spaces of $L$. Let  $AB$ be the $K$-vector space spanned by elements of the form $ab$ for $a\in A$ and $b \in B$.  Given $\dim_K A$ and $\dim_K B$, Hou-Leung-Xiang give a lower bound on  $\dim_K AB$, which we now describe.

Define $\Stab(AB)$ to be the set of elements in $L$ stabilizing $AB$ when acting by multiplication in the field:
$$
\Stab(AB) := \{ x \in L : x (AB) \subset AB \}.$$
Notice that $\Stab(AB)$ is an intermediate field extension of $L/K$:  it contains $K$, and is clearly  closed under field operations (where defined).

\tpoint{Lemma \cite[Thm.~2.4]{hlx}} {\em \label{l:hlx}We have
$$\dim_K AB \geq \dim_K A + \dim_K B - \dim_K \Stab(AB).$$}

\epoint{Proof of Theorem~\ref{t:one}}
Suppose otherwise that there are $i$ and $j$ such that $e_{i+j}> e_i + e_j$.  In Lemma~\ref{l:hlx}, let $K = K(\proj^1)$, let $L = K(C)$, let $A$ be the $K$-subvector space of $L$ generated by  $1$, $x_1$, \dots, $x_i$, and let $B$ be the $K$-subvector space generated by $1$, $x_1$, \dots, $x_j$.  From Observation~\ref{o:obs}, $AB$ is contained in the $K$-subvector space generated by $1$, $x_1$, \dots, $x_{i+j-1}$, so by Lemma~\ref{l:hlx} we have
$$
i+j \geq (i+1) + (j+1) - \dim_K \Stab (AB)$$
from which $\dim_K \Stab(AB) \geq 2$.  Because $\Stab(AB)$ is an intermediate field of $L/K$, strictly containing $K$, and by the primitivity hypothesis, there are no strictly intermediate fields, we must have $\Stab(AB)=L$, from which $AB=L$, contradicting the fact that $AB$ is contained in the vector subspace of $L$ generated by $1$, $x_1$, \dots, $x_{i+j-1}$. \epf

\section{The possible scrollar invariants do not form a convex set if $d$ is not a prime power}

\label{s:six}(This section does not rely on the results of \S \ref{s:inV}.  It is only placed here because of where it is in the narrative of the introduction.)

In this section, we show that if $d$ is not a prime power, then the possible scrollar invariants do not form a convex set for infinitely many $g$.  With some care this argument should be modifiable to show the result for all  $g \gg_d 1$ (for $d$ not a prime power), but it did not seem worth the trouble.  For small $g$, nonconvexity could plausibly ``accidentally'' fail because of the paucity of integer points.  The argument is a bit intricate, so for the sake of comprehensibility, we first give a concrete example for $d=6$, then give the case when $d$ is the product of two distinct primes, and finally  give the general argument.

\bpoint{The case $d=6$}
Our strategy is as follows.  We will construct two sextic covers  $\pi: C \rightarrow \proj^1$, both by genus $g=994$ curves, which  have scrollar invariants $(1, 166, 167, 332, 333)$ and $(1, 2, 331, 332, 333)$ respectively.    We then show that there can be no such cover with scrollar invariant $(1, 111, 222, 332, 333)$, which lies on the line segment joining the previous two points.

\epoint{The two covers}
Each cover will be constructed as a compositum of a hyperelliptic (degree 2) and a cubic (degree 3) cover.  In particular, these covers are not primitive, which is of course necessary given Theorem~\ref{t:one}.

Choose a general homogeneous polynomial  $f(x,y)$ of degree $r_2$ (to be named later, where $2| r_2$).  Consider the hyperelliptic cover given by
\begin{equation}
\label{eq:zf}
z^2 = f(x,y),
\end{equation}
say $\pi_2: C_2 \rightarrow \proj^1$.  Then $(\pi_2)_* \oh_{C_2} \cong \oh \oplus \oh(-r_2/2)$, and indeed we should interpret the $\oh_{\proj^1}$-algebra $(\pi_2)_* \oh_{C_2}$ as $\oh_{\proj^1} \oplus \oh_{\proj^1}(- r_2/2) z$, with multiplication induced by the relation \eqref{eq:zf}, i.e.,  multiplication in the $\oh_{\proj^1}$-algebra 
$$\left( \oh_{\proj^1} \oplus \oh_{\proj^1}(- r_2/2) z \right)
\times \left( \oh_{\proj^1} \oplus \oh_{\proj^1}(- r_2/2) z \right) \rightarrow \left( \oh_{\proj^1} \oplus \oh_{\proj^1}(- r_2/2) z \right)$$ is given by $z \times z \rightarrow (f(x,y), 0)$.

Similarly choose a general homogeneous polynomial  $g(x,y)$ of degree $r_3$ (to be named later, where $3|r_3$).  Consider the trigonal cover given by 
\begin{equation}
w^3=g(x,y),\label{eq:wg}
\end{equation}
say $\pi_3: C_3 \rightarrow \proj^1$, so we similarly have an isomorphism of $\oh_{\proj^1}$-algebras
$$
(\pi_3)_* \oh_{C_3} \overset \sim \longleftrightarrow \oh_{\proj^1} \oplus \oh_{\proj^1}(-r_3/3) w \oplus \oh_{\proj^1}(-2 r_3/3) w^2
$$
with multiplication induced by the relation \eqref{eq:wg}.

Now consider the sextic cover $\pi: C \rightarrow \proj^1$ given by both \eqref{eq:zf} and \eqref{eq:wg} (i.e., the compositum), so we have an isomorphism of  $\oh_{\proj^1}$-algebras
\begin{eqnarray*}
\pi_* \oh_{C} & \overset \sim \longleftrightarrow & \oh_{\proj^1} \oplus \oh_{\proj^1}(-r_3/3) w \oplus \oh_{\proj^1}(-2 r_3/3) w^2 \\
& & 
\oplus \oh_{\proj^1}(-r_2/2) z \oplus \oh_{\proj^1}(-r_2/2 -r_3/3) z w \oplus \oh_{\proj^1}(-r_2/2 -2 r_3/3) z w^2
\end{eqnarray*}
with multiplication induced by both relations \eqref{eq:zf} and \eqref{eq:wg}.  By the generality of our choice of $f(x,y)$ and $g(x,y)$, they have isolated roots distinct from each other, and it is straightforward to check that $C$ is smooth.    We thus visibly see the scrollar invariants over the cover $\pi$ as $(r_3/3, 2 r_3/3, r_2/2, r_2/2 +r_3/3, r_2/2 +2 r_3/3)$, ordered.    Choosing $r_2=2$ and $r_3=498$ (and reordering the summands)  yields the first cover, and $r_2=662$ and $r_3=3$ yields the second.

\epoint{The impossibility of  $(1, 111, 222, 332, 333)$}
Finally, we show that we cannot have a sextic cover $\pi: C \rightarrow \proj^1$ from an irreducible curve
with
\begin{equation}
\label{eq:no6}
\pi_* \oh_C \cong \oh_{\proj^1} \oplus \oh_{\proj^1}(-1) \oplus \oh_{\proj^1}(-111) \oplus \oh_{\proj^1}(-222) \oplus \oh_{\proj^1}(-332) \oplus \oh_{\proj^1}(-333).
\end{equation}
Consider the structure coefficients of the multiplication map of this $\oh_{\proj^1}$-algebra.  If $E_0= \oh$, \dots, $E_5= \oh(-e_5)$ are the summands of $\pi_* \oh_C$, then these coefficients are elements of $\Hom(E_i \otimes E_j, E_k) = H^0(\proj^1, \oh(e_i+e_j-e_k))$, which is necessarily zero if $e_i+e_j<e_k$.  

Now consider the field extension $K(C)$ over $K=K(\proj^1)$, with basis $1=x_0, x_1, \dots, x_5$ corresponding to the summands of \eqref{eq:no6}.  The structure coefficients are as described in the previous paragraph.  In particular, the $K$-vector space spanned by $x_0$ and $x_1$ is a subalgebra of $K(C)$ (as $1+1<111$, so $\Hom(E_i \otimes E_j, E_k) =0$ for $i,j \leq 1$, $k \geq 2$), and is thus necessarily a subfield $F$ of $K(C)$, of degree $2$ over $K$.  The $K$-vector space spanned by $x_0$, $x_1$, and $x_2$ is preserved by the multiplication by $F$ (as $1+111<222$, so $\Hom(E_i \otimes E_j, E_k)=0$ for $i \leq 1$, $j \leq 2$, $k \geq 3$), so it is a vector space over $F$.    But it is a 3-dimensional $K$-vector space, and thus cannot be a vector space over a degree $2$ field extension $F$ of $K$, contradiction.

\bpoint{The case $d=pq$ where $p<q$ are distinct primes}
We extend this idea to the case $d=pq$.  We first construct two covers, both as composita of cyclic covers of degree $p$ and $q$, and determine their scrollar invariants.  We then consider the ``midpoint'' of their scrollar invariants, and show that it cannot be achieved as the scrollar invariants of a cover of $\proj^1$ by an irreducible curve.

\epoint{The two covers} \label{s:twocov}Take $(a_1, b_1) = (2, 2+2N(p-1))$ and $(a_2, b_2) = (2+ 2N(q-1), 2)$ for $N \gg_{p, q} 1$.  
The $i$th  cover ($i=1,2$) is the compositum of a general cyclic cover of degree $p$ over a degree $a_i$ generally chosen branch locus and a general degree $q$ cyclic cover over a degree $b_i$ generally chosen branch locus. The scrollar invariants of the cyclic degree $p$ cover are $(a_i,\dots,(p-1)a_i)$ and the scrollar invariants of the cyclic degree $q$ cover are $(b_i,\dots,(q-1)b_i)$.  Let the two covers be denoted $\pi_1 \colon C_1 \rightarrow \PP^1$ and $\pi_2 \colon C_2 \rightarrow \PP^1$, and let $\pi \colon C \rightarrow \PP^1$ be the compositum.

The pushforward of the structure sheaf (with summands not necessarily written in increasing order) is
$$
(\pi_i)_* \oh_{C_i} \cong \bigoplus_{j=0}^{p-1} \bigoplus_{k=0}^{q-1}    \oh_{\proj^1}(- j a_i - k b_i)$$
which has degree
$(pq/2)  ( (p-1) a_i + (q-1) b_i)$, which is independent of $i$.

With summands written in decreasing order, the vector bundle $(\pi_1)_*\oh_{C_1}$ is:
\begin{equation}
\oh_{\proj^1} \oplus \oh_{\proj^1}(-2) \oplus \cdots \oplus \oh_{\proj^1}(-2(p-1)) \oplus \oh_{\proj^1} (-b_1) \oplus \oh_{\proj^1} (-2-b_1) \oplus \cdots.\label{eq:first}\end{equation}
Note that the summands at each step decrease by either $2$ or $b_1-2(p-1)$.
With summands written in decreasing order, the vector bundle $(\pi_2)_*\oh_{C_2}$ is:
\begin{equation}\label{eq:second}
\oh_{\proj^1} \oplus \oh_{\proj^1}(-2) \oplus \cdots \oplus \oh_{\proj^1}(-2(q-1)) \oplus  \oh_{\proj^1}(-a_2) \oplus \cdots .\end{equation}
The summands at each step decrease by either $2$ or $a_2-2(q-1)$.
The (vector bundle corresponding to the) average of the (lattice points corresponding to) the two scrollar invariants above begins as:
\begin{equation}\label{eq:third}
\oh_{\proj^1} \oplus \oh_{\proj^1}(-2) \oplus \cdots \oplus \oh_{\proj^1}(-2(p-1)) \oplus (-b_1/2-p) \oplus (-2-b_1/2-p) \oplus \cdots.\end{equation}
As in the special case $d=6$, the condition $N \gg 1$ implies that the first $p$ summands form a subalgebra ($2(p-1) + 2(p-1) < b_1/2+p$, so $\Hom(E_i \otimes E_j,E_k)=0$ for $i, j <p$ and $k \geq p$), so they correspond to a degree $p$ field extension $F$ of $K=K(\proj^1)$.   

\tpoint{Claim} {\em The first $q$ summands, pulled back to the generic point of $\proj^1$, are preserved by multiplication by elements of $F$.}

\bpf By our observations about how the degrees of the summands of \eqref{eq:first} and \eqref{eq:second} increase (written just after the two equation displays), in \eqref{eq:third} the summands satisfy \begin{eqnarray*}
-\deg E_q + \deg E_{q-1} &=& \frac {2 + (a_2 - 2(q-1))} 2
\\ &=& a_2/2 -q+2 \\ &=& 1 + N(q-1)-q+2 \\ &>& 2(p-1)  \quad \quad \text{(as $N \gg 1$)}\\ & =&  -\deg E_{p-1}.\end{eqnarray*}
Thus $\Hom(E_i \otimes E_j, E_k)=0$ for $i <p$, $j <q$, and $k \geq q$. \epf

Thus the first $q$ summands form a vector space over $F$.  But a $q$-dimensional vector space over $K$ cannot be a vector space over a degree $p$ field extension $F$ of $K$.

\bpoint{The general case $d=pqn$ where $p<q$ are distinct primes}  \label{s:pqn}The general case is a mild variation; we take composita of {\em three} cyclic field extensions, where the third is degree $n$ over a degree generally chosen branch locus with $M \gg_N 1$.  (We use the same $M$ in constructing both of the two covers analogous to \S \ref{s:twocov}.)  In the direct sum decomposition of the two covers constructed, the nontrivial summands of the third cyclic extension do not appear in the first $q$ summands, and thus do not interact with the previous argument.

\section{All concave scrollar invariants (a  positive proportion of the polytope  $\polytope_d$) are achieved}

\label{s:positive}
In this section, we prove Theorem~\ref{mainconstructionthm}. The degree of the cover is of course $d$, and the genus of $C$ is $\sum e_i - (d-1)$ (from \S \ref{s:early}). 

\point We begin by fixing some terminology.  Let $a$ be a nonnegative integer.  Suppose $C \subset \F_a$ is a curve in class $dD + kF \in \Pic \F_a$, where the base $\proj^1$ has projective coordinates $s$ and $t$.  Suppose $C$ has  equation (informally speaking, in ``toric coordinates'') \begin{equation}
\label{eq:f}f = f_0x^d + \dots + f_d y^d,
\end{equation}
where $f_i \in \C[s,t]$ is homogeneous of degree $k+ia$.  Here $y=0$ is the directrix of $\F_a$.  The following theorem is the key input that gives us direct control of the scrollar invariants of $C$.

\tpoint{Theorem} \label{t:nwood} {\em 
\label{t:nakagawawood} With the notation as above,  $\pi_* \oh_{C} \cong \oh_{\PP^1} \oplus \oh_{\PP^1}(-k-a) \oplus \oh_{\PP^1}(-k-2a) \oplus \dots \oplus \oh_{\PP^1}(-k-(d-1)a)$ with the following algebra structure.   Let $E_0,\dots,E_{d-1}$ be the $d$ summands.  For $i \leq j$, the multiplication map $E_i \otimes E_j \rightarrow E_k$ is given by:
\begin{eqnarray*}
  -f_{i + j - k} & & \text{if $\max\{i + j - d,1\} \leq k \leq i$} \\
f_{i + j - k} & & \text{if $j < k < d$} \\
   -f_{i + j - d}f_d & & \text{if $k = 0$ and $i + j \geq d$} \\
1  &  &\text{if $i = 0$ and $j = k$} \\
0 &  &\text{otherwise.}
\end{eqnarray*}}

This is a special case of \cite[\S 3]{wood}; see also \cite{na}. Now to prove Theorem~\ref{mainconstructionthm}, we describe a  particular cover $C' \rightarrow \PP^1$ whose normalization $C$ will have the desired Tschirnhausen bundle.  

We take $a=e_1$.  
Take a general curve $C'$ on $\F_a$ in class $dD$ (so $k=0$ in \eqref{eq:f}) subject to the constraint that $t^{ai-e_i}|f_i(s,t)$ ($0 \leq i \leq d$).  (Note that the assumptions of Theorem~\ref{mainconstructionthm}, namely concavity of the $e_i$, ensure that $ai-e_i \geq 0$ for all such $i$. Recall also that by convention, we set $e_0 = e_d = 0$.)

\tpoint{Claim}  {\em The curve $C'$ is irreducible, smooth away from $x=t=0$, and $C' \rightarrow \proj^1$ is primitive.} \label{cl:Cprimenice}

\bpf 
We consider the linear system of curves on the Hirzebruch surface 
with our imposed singularity at $x=t=0$ (a subsystem of $| \oh_{\F_a}(dD)|$).
Smoothness follows from Bertini's Theorem applied to the Hirzebruch surface minus the points $x=t=0$ (in the form of  \cite[Ex.~21.6.I]{risingsea}):  the linear system in question  on the Hirzebruch surface is base-point-free away from $x=t=0$.  It is easy to exhibit an irreducible curve in the system.  (We omit the details, because we will also soon show that if $C \rightarrow C'$ is the normalization, then $h^0(C, \oh_C)=1$, which also will give irreducibility.)  It is also straightforward to exhibit a curve in the linear system whose points above $[s;t] = [0;1] \in \PP^1$ are regular and branched with multiplicity $d-1$ and $1$, respectively (e.g., in local coordinates near $s=0$, $s + y^{d-1}(y-1)$); this implies primitivity.  \epf

Define the locally free sheaf  $\cA = \oplus_{i=0}^{d-1} t^{-ia+e_i} E_i$.  Then $\pi_* \oh_{C'} \subset \cA$ (as $e_i  \leq i e_1$ by concavity).  

\tpoint{Claim} \label{c:multiplicationextends} {\em The multiplication map of Theorem~\ref{t:nakagawawood} extends to give an $\oh_{\PP^1}$-algebra structure on $\cA$. } 

\bpf
For every $i$,$j$, and $k$, let $g_{ijk}$ represent the multiplication coefficient given in Theorem~\ref{t:nakagawawood}; it suffices to show that
\begin{equation}
\label{eqn:multiplication}
    t^{-ia + e_i}t^{-ja + e_j}g_{ijk} \in t^{-ka + e_k}\C[s,t].
\end{equation}

We check the five cases given in Theorem~\ref{t:nakagawawood}. The last two cases (when $i = 0$ and $j = k$, and when the multiplication coefficient is $0$) are trivial. We check the other three cases now. If $\max\{i+j-d,1\} \leq k \leq i$ or $j < k < d$ then $g_{ijk} = \pm f_{i+j-k}$. By comparing $t$-valuations, the inclusion (\ref{eqn:multiplication}) is equivalent to 
\[
    (-ia + e_i) + (-ja + e_j) + a(i + j - k) - e_{i + j - k} \geq -ka + e_k.
\]
The above is equivalent to $e_i + e_j \geq e_{i+j - k} + e_k$, which follows from concavity. If $k = 0$ and $i+j \geq d$, then $g_{ijk} =f_{i+j-d}f_d$. Proceeding as above, it suffices to show that
\[
    (-ia + e_i) + (-ja + e_j) + a(i + j - d) - e_{i + j - d} + ad - e_{d} + \geq 0.
\]
The above is equivalent to $e_i + e_j \geq  e_{i + j - d} + e_d$, which again follows from concavity. 
\epf

Thus
$$
C \coloneq \underline{\Spec} \;  \cA \rightarrow \underline{\Spec} \; \pi_* \oh_{C'} = C'
$$ 
is a birational morphism,  and an isomorphism away from $t=0$.

\tpoint{Claim} {\em The morphism $\underline{\Spec} \;  \cA \rightarrow \proj^1$ is unramified at $t=0$.}

\bpf
The ramification divisor of $C' \rightarrow \PP^1$ is given by the vanishing of $\Disc(f)$, and the ramification divisor of the map $C \rightarrow \PP^1$ is the previous quantity minus $2\sum_{i = 1}^{d-1} (ia - e_i)$ copies of $[1 \colon 0]$. Thus, to show that the ramification divisor of  $C \rightarrow \PP^1$ is not supported at $[1 \colon 0]$, it suffices to show that the power of $t$ dividing $\Disc(f)$ is bounded above by $2\sum_{i = 1}^{d-1} (ia - e_i)$.

Consider the universal binary $d$-ic form $F = F_0x^d + \dots + F_dy^d \in \ZZ[F_0,\dots,F_d]$. Lemma~\ref{cl:disc} shows that the polynomial $\Disc(F) \in \ZZ[F_0,\dots,F_d]$ contains the monomial $F_1^2\dots F_{d-1}^2$. Because our polynomial $f$ is general subject to the constraint that $t^{ia - e_i} \mid f_i$, this implies that the power of $t$ dividing $\Disc(f)$ is bounded above by $2\sum_{i = 1}^{d-1} (ia - e_i)$.
\epf 

\tpoint{Lemma} {\em
\label{cl:disc}
Consider the universal binary $d$-ic form $F = F_0x^d + \dots + F_dy^d \in \ZZ[F_0,\dots,F_d]$. The monomial $F_1^2\dots F_{d-1}^2$ occurs in the polynomial $\Disc(F)$.
}

\bpf
This is a completely formal statement.  The discriminant $\Disc(F)$ of
$F$ with respect to $(x,y)$  is a polynomial in the (symbolic) coefficients $F_i$
of total degree $2d-2$. The statement of the lemma is about a term in
this polynomial.  Let $G = F_1x^{d-2} + \dots + F_{d-1}y^{d-2}$.
 From \cite[p. 407, Equation 1.41]{gkz} we have:
\[
    \Disc(F)(0,F_1,\dots,F_{d-1},0) = \Disc(G)F_1^2F_{d-1}^2.
\]
The proof concludes via induction on $d$ (with base cases $d=2,3$ obvious).
\epf

Hence $C= \underline{\Spec} \; \cA$ is regular (hence normal) above $t=0$, so $C \rightarrow C'$ must be the normalization.  (Recall from Claim~\ref{cl:Cprimenice} that $C$ was normal away from $t=0$.)

\epoint{Completion of proof of Theorem~\ref{mainconstructionthm}}
To complete the proof of Theorem~\ref{mainconstructionthm}, we simply observe that $\cA = \pi_* \oh_C = \oplus_{i = 0}^{d-1} E_i$, where $E_i \cong \oh(-e_i)$. \epf

(We also now immediately see that $h^0(C, \oh_C) = \sum_{i = 0}^{d-1} h^0(\proj^1, \oplus \oh(-e_i))  = 1$, fulfilling a promise in the proof of Claim~\ref{cl:Cprimenice}.)

\bpoint{Generalizations of the construction} 
\label{r:variant}
In the proof of Theorem~\ref{mainconstructionthm}, we started with a concave sequence $(0,e_1,\dots,e_{d-1},0)$, considered certain binary forms $f$ such that $t^{ai-e_i} \mid f_i$, and showed that the scrollar invariants of the normalization were given by $(e_1,\dots,e_{d-1})$. More generally, given any $(a_0,\dots,a_d) \in \ZZ_{\geq 0}^{d+1}$ and a form $f$ such that $t^{a_i} \mid f_i$, we compute the scrollar invariants of the normalization below.

Let $f \in \Sym^d (\oh \oplus \oh(a)) \otimes \oh(k)$ and write 
\[
    f = f_0x^d + \dots + f_d y^d,
\]
again in ``toric coordinates'' where $f_i \in \CC[s,t]$ is a homogeneous polynomial of degree $k + ia$. The construction of \cite[Thm.~1.2]{wood} associates to $f$ a finite flat morphism $\pi: C' \rightarrow \proj^1$ (note that $f$ can even be the zero form). Let  $\{(i,b_i)\}_{0 \leq i \leq d}$ be the lower convex hull of $\{(i,a_i)\}_{0 \leq i \leq d}$. 

\tpoint{Theorem (``general $(a_0,\dots,a_d)$'')} {\em 
\label{p:generalization}
If $t^{a_i} \mid f_i$, then there exists a curve $C$ and a birational morphism $C \rightarrow C'$ such that the composition $C \rightarrow C' \rightarrow \PP^1$ has scrollar invariants
\[
    (\deg(f_1) - \lfloor b_1 \rfloor,\dots,\deg(f_{d-1}) - \lfloor b_{d-1} \rfloor).
\]
Now suppose $f$ is general subject to the constraint that $t^{a_i} \mid f_i$. If $\deg(f_i) \geq a_i$ for all $0 \leq i \leq d$ and $\deg(f_i) - b_i \geq 1$ for all $1 \leq i \leq d-1$, then $C$ is smooth and irreducible. 
}

\epoint{Remark: how we used Theorem~\ref{p:generalization} to prove Theorem~\ref{mainconstructionthm}} {
Theorem~\ref{mainconstructionthm} is simply the case of Theorem~\ref{p:generalization} where $k = 0$ and $a = e_1$, and $a_i = ai - e_i$. The concavity of $(0,e_1,\dots,e_{d-1},0)$ now implies that $b_i = ia - e_i$, so
\[
    \deg(f_i) - b_i = (k + ai) - (ia - e_i) = e_i.
\]
}

A natural question is: can we use Theorem~\ref{p:generalization} to produce even more scrollar invariants? The answer is yes: when $d = 6$, the scrollar invariant $(1,1,1,2,2)$ does \emph{not} arise from Theorem~\ref{t:one} as it is not concave (and no reordering is concave). However, it does arise from Theorem~\ref{p:generalization}: simply take $k = 2$, $a = 0$, and $(a_0,\dots,a_6) = (2,2,1,1,0,0,0)$. One can easily generalize this example to any degree $d \geq 4$. (There is also an example for $d=5$, but it is more complicated to explain.)

Sadly, although Theorem~\ref{p:generalization} allows us to construct scrollar invariants not present in Theorem~\ref{t:one}, it does not obtain a larger proportion of $\polytope_d$, as $g \rightarrow \infty$, as compared to the current approach. This can be checked via a lengthy computation, but this argument is not relevant to our discussion here. 

\epoint{Proof of Theorem~\ref{p:generalization}} 
The $\oh_{\PP^1}$-algebra structure of $\pi_* \oh_{C'}$ is described in Theorem~\ref{t:nwood} (see \cite[\S 3]{wood}). Define the locally free sheaf
\[
    \cA = \oh_{\PP^1} \oplus t^{-\lfloor b_1 \rfloor}\oh_{\PP^1}(-e_1) \oplus \dots \oplus t^{-\lfloor b_{d-1} \rfloor}\oh_{\PP^1}(-e_{d-1}).
\]
Observe that $\pi_* \oh_{C'} \subseteq \cA$. The multiplication map on $\pi_* \oh_{C'}$ extends to an $\oh_{\PP^1}$-algebra structure on $\cA$; this can be seen by easily generalizing the computation in \S \ref{c:multiplicationextends}. Thus as before, 
$$
C \coloneq \underline{\Spec} \;  \cA \rightarrow \underline{\Spec} \; \pi_* \oh_{C'} = C'
$$ 
is a birational morphism,  and an isomorphism away from $t=0$. This completes the proof of the first part of the theorem. 

Now suppose $\deg(f_i) \geq a_i$ for all $0 \leq i \leq d$ and $f$ is general subject to the constraint that $t^{a_i} \mid f_i$. The resulting linear system is base-point-free away from $t = 0$, and thus $C$ is smooth away from $t = 0$. Consider the commutative diagram:
\[
\begin{tikzcd}
C_\ell \arrow{r}{\gamma} \arrow[swap]{d}{\rho} & C \arrow{d}{\phi}\\
C'_\ell \arrow{r}{} \arrow[swap]{d}{\tau} & C' \arrow{d}{\pi} \\
\PP^1 \arrow{r}{m_\ell} & \PP^1
\end{tikzcd}
\]
where $\ell$ is the least common multiple of the denominators of $b_i$, and $m_\ell \colon \PP^1 \rightarrow \PP^1$ sends $[u : v] \rightarrow [u^\ell : v^\ell]$, and $C'_\ell$ is the pullback of the bottom square. Let $C_\ell \rightarrow C'_\ell$ and $C \rightarrow C'$ be the partial normalization maps as described above. Note that the top square need not be a pullback diagram. 

\tpoint{Claim} {\em The morphism $\tau \circ \rho \colon C_{\ell} \rightarrow \proj^1$ is unramified above $u = 0$.}

\bpf
Because $C'$ is defined by an equation of the form
\[
    f(x,y) = t^{a_0}g_0(s,t)x^d + \dots + t^{a_d}g_d(s,t)y^d,
\]
we have that $C'_\ell$ is defined by an equation of the form
\[
    f_\ell(x,y) = u^{\ell a_0}g_0(u^\ell,v^\ell)x^d + \dots + u^{
    \ell a_d}g_d(u^\ell,v^\ell)y^d
\]

Let $F = U^{\ell a_0}G_0x^d + \dots + U^{\ell a_d}G_dy^d \in \CC[U,G_0,\dots,G_d]$ be the universal such form.

\tpoint{Sub-claim} { \em We have $\Disc(F) = U^{2\sum_{i = 1}^{d-1}\ell b_i}H(U,G_0,\dots,G_d)$
for a polynomial $H \in \CC[U,G_0,\dots,G_d]$ which is not divisible by $U$. 
}

\bpf
Let $v_U$ be the $U$-adic valuation on $\CC[U,G_0,\dots,G_d]$. It suffices to show that
\[
    v_U(\Disc(F)) = 2\sum_{i = 1}^{d-1}\ell b_i.
\]
Let $\alpha_1,\dots,\alpha_d$ be the roots of $F(x,1)$, ordered so that $v_U(\alpha_1)\leq \dots \leq v_U(\alpha_d)$. We have that 
\[
    \Disc(F) = \left(U^{\ell a_0}\right)^{2d - 2} \prod_{i < j}(\alpha_i - \alpha_j)^2.
\]
Because the lower convex hull of the points $(i,\ell a_i)$ is exactly $(i,\ell b_i)$, we have $a_0 = b_0$. Moreover, by the theory of Newton polygons, we have $v_U(\alpha_i) = \ell(b_i - b_{i - 1})$. Moreover, because $F$ is the universal such form, we have
\begin{equation}
\label{eqn:valuations}
 v_U(\alpha_i - \alpha_j) = \min(v_U(\alpha_i),v_U(\alpha_j)).
\end{equation}
Computing $U$-adic valuations, we see that
\begin{equation}
\label{eqn:valuations-2}
v_U(\Disc(F)) = (\ell a_0)(2d - 2) + 2\sum_{i < j}\min(v_U(\alpha_i),v_U(\alpha_j)).
\end{equation}
Expanding out the last term, we see that:
\begin{equation}
\label{eqn:valuations-3}
\sum_{i < j}\min(v_U(\alpha_i),v_U(\alpha_j)) = \sum_{i < j} v_U(\alpha_i) =\ell \sum_{i < j} (b_i - b_{i - 1}) = \ell\sum_{i = 1}^{d-1}(d-i)(b_i - b_{i - 1}).
\end{equation}
Now, a short computation shows that
\[
    (\ell b_0)(2d - 2) + 2\ell\sum_{i = 1}^{d-1}(d-i)(b_i - b_{i - 1}) = 2\sum_{i = 1}^{d-1}\ell b_i
\]
as required.
\epf

Evaluating on $f_{\ell}$ we obtain: 
\[
    \Disc(f_{\ell}) = u^{2\sum_{i = 1}^{d-1}\ell b_i}H(u,g_0(u^\ell,v^{\ell}),\dots,g_d(u^\ell,v^{\ell})).
\]
Because $\deg(f_i) \geq a_i$, all the $g_i$ have nonnegative degree, and so for a general choice of $g_0,\dots,g_d$, the polynomial $H(u,g_0(u^\ell,v^{\ell}),\dots,g_d(u^\ell,v^{\ell}))$ will not be divisible by $u$. Therefore the ramification divisor of $\tau : C_{\ell}' \rightarrow\PP^1$ has degree precisely $2\sum_{i = 1}^{d-1}\ell b_i$ at $u = 0$. 

The scrollar invariants of $\tau : C_{\ell}' \rightarrow\PP^1$, again by Theorem~\ref{t:nwood}, are $(\ell\deg(f_1),\dots,\ell\deg(f_{d-1}))$. Now, by definition the scrollar invariants of $\tau \circ \rho \colon C_{\ell} \rightarrow \PP^1$ are:
\[
    (\ell\deg(f_1) - \ell b_1,\dots,\ell\deg(f_{d-1}) - \ell b_{d-1}).
\]
Hence the degree of the ramification divisor of $\tau \circ \rho \colon C_{\ell} \rightarrow \PP^1$ at $u = 0$ is $2\sum_{i = 1}^{d-1}\ell b_i$ less than that of $C_{\ell}'$, and hence $C_{\ell}$ is unramified above $u = 0$. 
\epf

\bpoint{Completion of the proof of Theorem~\ref{p:generalization}}
Because $\tau \circ \rho \colon C_{\ell} \rightarrow \proj^1$ is unramified, the points of $C_{\ell}$ above $t = 0$ in the map $m_{\ell} \circ \tau \circ \rho = \pi \circ \phi \circ \gamma$ are smooth. 

There is a natural inclusion $\oh_C \subseteq \gamma_* \oh_{C_\ell}$ which pushes forward to an inclusion
\[
    (\pi \circ \phi)_*\oh_C \subset (\pi \circ \phi \circ \gamma)_* \oh_{C_\ell}
\]
Restricting to the affine line $\AAA^1 = \proj^1 \setminus \infty$ and taking global sections, we have that
\[
     \Gamma(\AAA^1, (\pi \circ \phi)_*\oh_C) \subset \Gamma(\AAA^1, (\pi \circ \phi \circ \gamma)_* \oh_{C_\ell}).
\]
By construction,
\[
    \Gamma(\AAA^1, (\pi \circ \phi)_*\oh_C) = \Gamma(\AAA^1, (\pi \circ \phi \circ \gamma)_* \oh_{C_\ell}) \cap K(C).
\]
The last equality shows that the smoothness of $C_{\ell}$ implies the smoothness of $C$ above $t = 0$. 

Finally, $h^0(C, \oh_C) = h^0(\proj^1, \oh) + \sum_{i = 1}^{d-1} h^0(\proj^1, \oplus \oh(b_i - \deg(f_i)))$. If $\deg(f_i) - b_i \geq 1$ for all $1 \leq i \leq d - 1$, then $h^0(C, \oh_C) = 1$ so $C$ is irreducible.
\epf

There is one more trick we may play to construct \emph{even more} scrollar invariants --- we may start with singularities above multiple points of $\PP^1$. More precisely if $\ell_1,\dots,\ell_m \in \CC[s,t]$ are linear forms with distinct roots, then in the notation as in Theorem~\ref{p:generalization}, we obtain the following theorem, whose proof is the obvious generalization of that of Theorem~\ref{p:generalization}.

\tpoint{Theorem (``multiple points of the base'')} \label{t:furthergeneralization} { \em
Fix $(a^1_0,\dots,a^1_d),\dots,(a^m_0,\dots,a^m_d) \in \ZZ_{\geq 0}^{d+1}$ and let $b_i^j$ be the lower convex hull. If $\ell_j^{a_i^j} \mid f_i$ for all $i,j$, then there exists a curve $C$ and a birational morphism $C \rightarrow C'$ such that the composition $C \rightarrow C' \rightarrow \PP^1$ has scrollar invariants
\[
    (\deg(f_1) - \sum_{j} b_1^j,\dots,\deg(f_{d-1}) - \sum_{j} b_{d-1}^j).
\]
Now suppose the $\ell_j$ are general and $f$ is general subject to the constraint that $\ell_j^{a^j_i} \mid f_i$. If $\deg(f_i) \geq \sum_{j}a^j_i$ for all $0 \leq i \leq d$ and $\deg(f_i) - \sum_j b^j_i \geq 1$ for all $1 \leq i \leq d-1$, then $C$ is smooth and irreducible. 
}

\epoint{Proof of Theorem~\ref{t:furthergeneralization}} 
The proof of Theorem~\ref{t:furthergeneralization} is almost identical to that of Theorem~\ref{p:generalization}. The only additional difficulty is checking smoothness above the vanishing of the $\ell_j$; in other words, we need to prove that if the $\ell_j$ are general and $f$ is general subject to the given constraint and $\deg(f_i) \geq \sum_{j}a^j_i$ for all $0 \leq i \leq d$ and $\deg(f_i) - \sum_j b^j_i \geq 1$ for all $1 \leq i \leq d-1$, then $C$ is smooth. Suppose without loss of generality that $\ell_1 = t$. 

Following the idea behind the proof of Theorem~\ref{p:generalization}, let $\ell$ be the least common multiple of the $b^1_i$. Consider the commutative diagram as in the proof of Theorem~\ref{p:generalization} where now $C_{\ell}$ is the given partial normalization above $\ell_1$. Then $C_{\ell}'$ is defined by the equation
\[
    f_{\ell}(x,y) = u^{\ell a_0^1}\Big(\prod_{j = 2}^m\ell_j(u^\ell,v^{\ell})^{a_0^j}\Big)g_0(u^\ell,v^{\ell})x^d + \dots + u^{\ell a_d^1}\Big(\prod_{j = 2}^m\ell_j(u^\ell,v^{\ell})^{a_d^j}\Big)g_d(u^\ell,v^{\ell})y^d.
\]

As before, let $F \in \CC[U,L_2,\dots,L_m,G_0,\dots,G_d]$ be the universal such form. We have that
\[
    \Disc(F) = U^{2\sum_{i = 1}^{d-1}\ell b_i^1}H(U,L_2,\dots,L_m,G_0,\dots,G_d)
\]
for some polynomial $H$ which is not divisible by $U$. Evaluating on $f_{\ell}$, we get
\[
    \Disc(f_{\ell}) = u^{2\sum_{i = 1}^{d-1}\ell b_i^1}H(u,\ell_2(u^\ell,v^{\ell}),\dots,\ell_m(u^\ell,v^{\ell}),g_0(u^{\ell},v^\ell),\dots,g_d(u^\ell,v^\ell)).
\]
Because $\deg(f_i) \geq \sum_{j}a^j_i$ for all $0 \leq i \leq d$, we have that $\deg(g_i) \geq 0$. Similarly, $\deg(\ell_j) = 1$, so a general choice of $\ell_j$ and $g_i$ ensures that the polynomial 
\[
H(u,\ell_2(u^\ell,v^{\ell}),\dots,\ell_m(u^\ell,v^{\ell}),g_0(u^{\ell},v^\ell),\dots,g_d(u^\ell,v^\ell))
\]
is indivisible by $u$. Proceeding as in Theorem~\ref{p:generalization} completes the proof of the theorem.
\epf

Remarkably, Theorem~\ref{t:furthergeneralization} \emph{does} produce a new positive proportion of $\polytope_d$, but we are still unable to fill  $\polytope_d$ for all $d \geq 5$. For example, suppose we wanted to obtain the scrollar invariants $(3,3,3,5,5)$. It is \emph{not} possible to get these scrollar invariants by using only Theorem~\ref{p:generalization}. Indeed if we apply Theorem~\ref{p:generalization} to $k = 5$, $a = 0$, and $(a_0,\dots,a_6) = (0,0,0,2,2,4,4)$, and then the corresponding normalization will have the scrollar invariants $(2,3,4,5,5)$. However if we instead apply Theorem~\ref{t:furthergeneralization} with $m = 2$ and $(a^j_0,\dots,a^j_6) = (0,0,0,1,1,2,2)$ for $j= 1,2$, then we obtain the desired scrollar invariants. 

\section{The image of $\Tsch: \Hur_{d,g} \rightarrow \Bun_{\proj^1}$ is not preserved under generization ($d > 4, g \gg_d 1$)}
\label{s:notflat}

In low genus ($2g+1\leq d$), the image of $\Tsch: \Hur_{d,g} \rightarrow \Bun_{\proj^1}$ is preserved under generization. This is because the image of the map $\Tsch$ is a single point, and the corresponding bundle $\oh(e_1) \oplus \dots \oplus \oh(e_{d-1})$ has no nontrivial generization because $e_{d-1} - e_1 \leq 1$ (see Proposition~\ref{p:triangle}). Likewise, in low degree ($d \leq 4$), the image of $\Tsch$ is also preserved under generization. This is trivial for $d = 2$, and follows from \S \ref{ss:cubic} and Theorem~\ref{t:tet} respectively when $d = 3,4$. In this section we prove that when $d \geq 5$ and $g \gg_d 1$, the image of $\Tsch$ is \emph{not} preserved under generization, i.e. that there exists a smooth curve whose Tschirnhausen bundle admits a generization which does \emph{not} arise as the Tschirnhausen bundle of a smooth curve.  Along the way, we also prove that when $d \geq 4$ and $g \gg_d 1$, the image of $\Tsch$ restricted to primitive covers is also not preserved under generization.

\tpoint{Proposition} \label{p:notclosed}  {\em If $d \geq 5$, then for $g \gg_d 1$, the image of $\Tsch: \Hur_{d,g} \rightarrow \Bun_{\proj^1}$ is not preserved under generization.}

\bpf
Suppose $(e_1, e_2, \dots, e_{d-1})$ is a concave nondecreasing sequence of positive integers, with $e_1 = a$, $e_2 = 2a$, $e_3 = 3a$, $e_4 \geq 3a+5$ (where $a \in \Z$, $a\geq 5$). Then by Theorem~\ref{mainconstructionthm}, there is a primitive cover $\pi: C \rightarrow \proj^1$ by a smooth curve with Tschirnhausen bundle $\cE_\pi \cong \oplus_{i=1}^{d-1} \oh(e_i)$. If $(e'_1, \dots, e'_{d-1}) = (e_1, e_2+1, e_3+2, e_4-3, e_5, \dots, e_{d-1})$, then  the bundle $\oplus_{i=1}^{d-1} \oh(e'_i)$ is {\em not} the Tschirnhausen bundle of a cover (Lemma~\ref{l:inequalities}), but it is a generization of $\cE_{\pi}$. The statement now follows from Proposition~\ref{p:fun}.
\epf

\tpoint{Lemma} \label{l:inequalities} { Let $C \rightarrow \PP^1$ be a cover of degree $\geq 4$ with scrollar invariants $(e_1,\dots,e_{d-1})$. If $e_2 > 2e_1$, then $e_3 \leq e_1 + e_2$.
}

\bpf
As in \S \ref{point:coords}, choose a coordinate $t$ on $\AAA^1  = \proj^1 \setminus\{\infty\}$ and a splitting 
\[
    \pi_* \oh_C \simeq \oh \oplus \oh(-e_1) \oplus \dots \oplus \oh(-e_{d-1}).
\]
Let $K = K(\proj^1)$ and let $L = K(C)$. As in \ref{point:coords}, this splitting gives rise to a basis $1,x_1,\dots,x_{d-1}$ for the vector space $L/K$.

Let $A$ be the $K$-subvector space spanned by $1,x_1$. Because $e_2 > 2e_1$, we have that $AA \subseteq A$, so $A$ is a quadratic field extension of $K$. Now let $B$ be the $K$-subvector space spanned by $1,x_1,x_2$. Because $AB$ is an $A$-vector space and $\dim_K A = 2$, the quantity $\dim_K AB$ must be even. Therefore, $\dim_K AB \geq 4$. Consequently there exists some integer $k \geq 3$ such that $e_k \leq e_1 + e_2$.
\epf

\tpoint{Proposition} \label{p:fun}  {\em Fix an integer $d \geq 5$.  Let $S$ be the set of length $d-1$ nondecreasing concave sequences $(e_1, e_2, \dots, e_{d-1})$ starting with $e_1 = a$, $e_2 = 2a$, $e_3 = 3a$, $e_4 \geq 3a+5$ (where $a \in \Z$, $a\geq 5$).  Let $T = \{ \sum e_i : (e_i) \in S \}$.   Then $\Z^+ \setminus T$ is a finite set.}

\bpf
For any $a \geq 2$, and $(d-2)a+5 \leq e_{d-1} \leq (d-1)a$, the sequence $(e_1=a, e_2=2a, \dots, e_{d-2} = (d-2)a, e_{d-1})$ is in $S$.  Thus for all $a \geq 5$,  the interval $[f(a), g(a)] \subset \Z^+$ is in $T$ for $f(a) = a+ \cdots + (d-2)a + ( (d-2)a+5) = \binom {d-1} 2 a + ((d-2)a + 5)$ and $g(a) = a+ \cdots + (d-2)a + (d-1)a = \binom {d-1} 2 a + (d-1) a$.  But it is easy to check that (i) $f(a+1)<g(a)$ for $a \gg_d 1$ (as $g(a)-f(a+1) = a -h(d)$ for appropriate $h(d)$), and (ii) $f(a) \rightarrow \infty$ as $a \rightarrow \infty$, so these intervals cover all but finitely many of the positive integers. \epf

\epoint{Remark: the set of Tschirnhausen bundles of primitive covers is not preserved under generization ($d \geq 4$ and $g \gg_d 1$)} 
\label{r:deformingprimitives}
Let $\Hur_{d,g}^p$ denote the open substack of primitive covers, and let $\Tsch^p: \Hur^p_{d,g} \rightarrow \Bun_{\proj^1}$ denote the restriction. We now show that when $d \geq 4$ and $g \gg_d 1$, the image of $\Tsch^p$ is not preserved under generization. 

Indeed, when $d = 4$, take any nondecreasing sequence $(e_1,e_2,e_3)$ with $e_1 = a$, $e_2 = 2a$, and $a+2 \leq e_3 \leq 3a$ for $a \geq 1$. Because $e_3 \leq e_1 + e_2$ and $e_2 \leq 2e_1$, it can be realized as the scrollar invariants of a smooth primitive cover (Theorem~\ref{t:tet}). Let $(e'_1,e'_2,e'_3) = (e_1,e_2 + 1,e_3-1)$ and note that the corresponding bundle does \emph{not} arise as the Tschirnhausen bundle of a primitive quartic cover (Theorem \ref{t:tet}), but is indeed a generization of $\oplus_{i = 1}^3 \oh(e_i)$. To conclude in this case it suffices to show that if $S$ is the set of length $3$ sequences $(e_1,e_2,e_3)$ as above, then the set $\ZZ^+\setminus \{\sum e_i : (e_i) \in S\}$ is finite. When $d \geq 5$, the proof of Proposition~\ref{p:notclosed} already shows that there exists a smooth \emph{primitive} cover whose Tschirnhausen bundle has a generization which is not the Tschirnhausen bundle of any smooth cover.

\section{Motivation and analogy from number theory}
\label{s:density}

\label{s:arithmetic}There is a particularly strong analogy to this question in number theory.   The analogy to $\PP^1$ is 
$\operatorname{Spec} \ZZ$, and the analogy to the curve $C$ is the ring of integers $\oh_K$ in a degree $d$ extension $K$ of $\QQ$; see \cite[\S 7]{hess} for more details on the arithmetic side.   

Our arguments in this paper extend readily to this arithmetic setting, but do not yield the same strength of results (number fields whose "shape" is in this "concave region"), because of one missing step:  the need for a sieving result to ensure that there is a maximal order in a number field, which is strictly harder than the geometric analogue of Bertini's Theorem (a blunt dimensional statement).

\bpoint{Extension:  Density functions on the polytopes}

We now make precise the connection we predict between the geometry of degree $d$ number fields, and the geometry of degree $d$ covers of $\proj^1$. Fix a positive integer $d \geq 2$.

\epoint{Definition} {
\label{def:pigeo}
For each open ball $B \subset \R^{d-1}$ and $g>0$, let 
$$\pi^{geo}(B,g) = \frac{1}{\dim \mathcal{H}_{d,g}}\dim  \left ( \bigcup_{\vec{e} \in \Z^{d-1} : \vec{e}/(d+g-1) \in B}    \mathcal{H}_{\vec{e},g} \right )$$
where $\mathcal{H}_{\vec{e}, g}$ is the moduli space of genus $g$ degree $d$ covers of $\proj^1$  with scrollar invariants $\vec{e}$. (The dimension of a finite disjoint union is  the maximum of the dimensions, and the dimension of the empty set is $-1$.) 
}

We now define the arithmetic analogue of $\pi^{geo}(B,g)$. For a degree $d$ number field $K$, let $\lambda_0,\dots,\lambda_{d-1}$ be the successive minima with respect to the following quadratic form: embed $K \xhookrightarrow{} \CC^d$ and take $\frac{1}{d}$ times the standard Hermitian quadratic form on $\CC^d$. (Or even better: consider $\oh_K$ as a metrized line bundle on $\Spec(\oh_K)$. Pushing forward via the degree $d$ cover $\Spec(\oh_K) \rightarrow \Spec(\ZZ)$ gives a metrized vector bundle on $\Spec(\ZZ)$. Then $\log \lambda_0,\dots,\log \lambda_{d-1}$ are simply the slopes of this metrized vector bundle.)

Then $\lambda_0 = 1$ and $\prod_{i = 1}^{d - 1}\lambda_{i} \asymp_d \sqrt{\lvert \Disc(K) \rvert}$ by the AM-GM inequality and Minkowski's second theorem respectively. The quantity $\log \lvert \Disc K \rvert$ is analogous to $2g$ and the quantities $\log \lambda_i$ are analogous to the scrollar invariants $e_i$. 

\epoint{Definition} {
For each open ball $B \subset \R^{d-1}$ and $X>0$, let 
\[
    \pi^{arith}(B,X) = \frac{\log \#\{K \; \mid \; \lvert \Disc(K) \rvert \in [X,2X], \;
    \log_{\sqrt{X}}(\lambda_1,\dots,\lambda_{d-1}) \in B\}}{\log\#\{K \; \mid \; \lvert \Disc(K) \rvert \in [X,2X]\}}
\]
(By definition, we take $\log(0) = -1$, just as the dimension of the empty set is taken to be $-1$, cf. Definition~\ref{def:pigeo}.)
}

Let $\Xi = \{ (\overline{e}_1, \dots, \overline{e}_{d-1}) \in \R^{d-1} : \sum \overline{e}_i = 1 \}$ be the hyperplane whose coordinates sum to $1$.

\tpoint{Conjecture}
\label{conj:measure}
{\em For all $B$, both $\lim_{g \rightarrow \infty} \pi^{geo}(B,g)$ and $\lim_{X \rightarrow \infty} \pi^{arith}(B,X)$ exist, are equal, and depend only on the intersection of $B \cap \Xi$ and not on $B$. 
}

Conjecture~\ref{conj:measure} is trivially true in the case $d = 2$. We give a proof in the case $d = 3$ below. 

\begin{figure}[ht]
\centering
\includegraphics[scale=0.35]{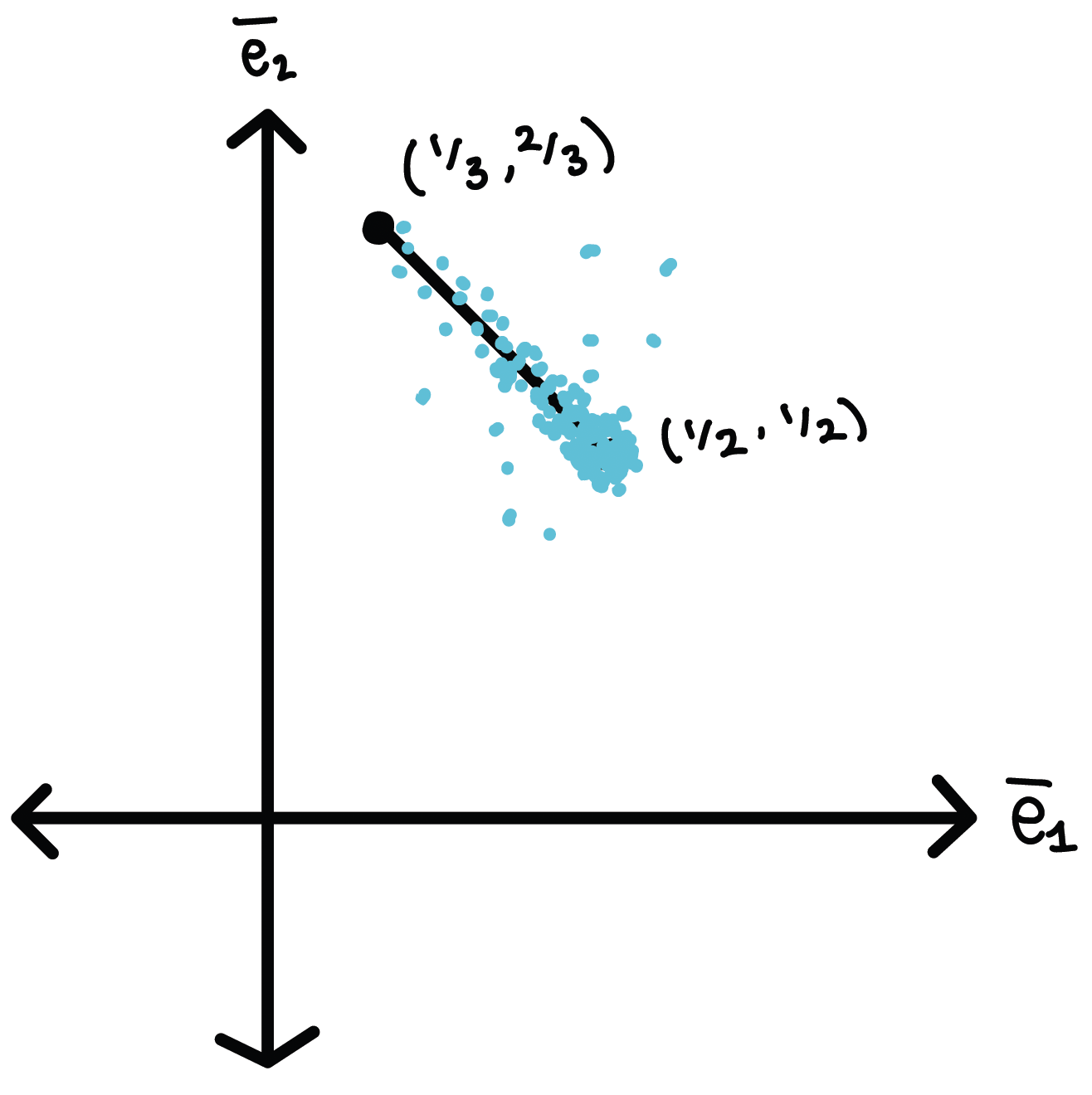}
\caption{An artist's depiction of the points $(\log_{\sqrt{X}}\lambda_1,\log_{\sqrt{X}}\lambda_2)$ as we range over cubic fields with absolute discriminant in the interval $[ X, 2X]$. Note how the points cluster close to the line segment, but are not necessarily contained on it. Moreover, the points are concentrated near the point $(1/2,1/2)$, and become more sparse as we move to the left hand side of the line segment. Beyond the left hand side of the line segment, the points abruptly end.}\label{f:cubicoffline}\end{figure}

\tpoint{Proposition}
\label{p:cubicconj}
{\em Conjecture~\ref{conj:measure} is true when $d = 3$.}

\bpf 
The locus $\cH_{\vec{e},g}$ is nonempty precisely when $e_1 + e_2 = g+2$, $e_1 \leq e_2$, and $e_2 \leq 2e_1$. In this case, the dimension of $\cH_{\vec{e},g}$ is given by
\[
    h^0(\proj^1, \Sym^3(\oh(\vec{e})) \otimes \oh(\vec{e})^{\vee}) - h^0(\proj^1,\End(\oh(\vec{e}))). 
\]
Moreover, the dimension of $\cH_{3,g} = 2g + 5$. So, if $B$ nontrivially intersects the closed line segment from $(1/3,2/3)$ to $(1/2,1/2)$, then
\[
    \lim_{g \rightarrow \infty}\pi^{geo}(B,g) = \sup \left\{1 - \frac{\ove_2-\ove_1}{2} : (\ove_1,\ove_2) \in B \cap \Xi \right\}
\]
If $B \cap \Xi = \emptyset$, then $\lim_{g \rightarrow \infty}\pi^{geo}(B,g) = 0$. The equivalent computation in the arithmetic case can essentially be found in the proof of \cite[Thm.~1.6]{sameera3}. 
\epf

For $d = 4$, the difficulty of sieving out maximal orders is already too hard in the arithmetic setting. For case of not necessarily maximal orders when $d = 4,5$, see \cite{sameera3} for a description of $\pi^{arith}$. In these cases, the nonmaximal version of $\pi^{arith}$ is described as the maximum of a continuous piecewise linear function supported on a polytope; in these known cases, this matches the version of $\pi^{geo}(B,g)$ including singular curves.

Let $\pi^{geo}(B)$ and $\pi^{arith}(B)$ be the corresponding limits. One might be tempted to ask: how does the density of curves or number fields concentrate on specific points of $\Xi$?  The following conjectural density function is meant to capture that information. For a point $x \in \Xi$ define
\[
    \rho^{geo}(x) = \lim_{i \rightarrow \infty}\pi^{geo}(B_i(x))
\]
where $\{B_i(x)\}$ is a set of open balls all containing $x$ and with radius converging to $0$. Define $\rho^{arith}$ similarly.

\tpoint{Conjecture}
\label{conj:rho}
{\em The quantities $\rho^{geo}(x)$ and $\rho^{arith}(x)$ exist and are equal. They are supported on a finite union of rational polytopes, are piecewise linear, and are continuous on their support.
}

\tpoint{Corollary} {\em
If $d = 3,4,5$, Conjecture~\ref{conj:rho} is true. 
}

\bpf
The case $d = 3$ follows from the explicit description of $\pi^{geo}(B_i(x))$ given in Proposition~\ref{p:cubicconj}. The case $d = 4$ follows from Theorem~\ref{t:quartic} (see Proposition~\ref{p:quartic-rho}). The case $d = 5$ can similarly be extracted from \cite[Theorem~1.5]{fv}.
\epf

\bpoint{Further extension: refined density functions taking into account ``all'' scrollar invariants}
As we saw in Section~\ref{s:quartic}, a quartic cover has three scrollar invariants $(e_1,e_2,e_3)$ and an associated trigonal cover with scrollar invariants $(f_1,f_2)$, with the property that $e_1 + e_2 + e_3 = f_1 + f_2 = g+3$. In this case, we can consider \emph{refined} density functions which take into account the extra data of $(f_1,f_2)$.

\epoint{Definition} { \em
For each open ball $B \subset \R^{5}$ and $g>0$, let
$$\pi_{ref}^{geo}(B,g) = \frac{1}{\dim \mathcal{H}_{4,g}}\dim  \left ( \bigcup_{(\vec{e},\vec{f}) \in \Z^{5} : (\vec{e},\vec{f})/(g+3) \in B}    \mathcal{H}_{\vec{e},\vec{f}} \right )$$
where $\cH_{\vec{e},\vec{f}} \subseteq \cH_{4,g}$ is the moduli space of genus $g$ degree $d$ covers of $\proj^1$  with scrollar invariants $\vec{e}$ and $\vec{f}$. (As before, the dimension of a finite disjoint union is  the maximum of the dimensions, and the dimension of the empty set is $-1$.) 
}

Similarly, one can define $\pi_{ref}^{arith}$ (in this case, we replace the data of $(f_1,f_2)$ with the successive minima of the cubic resolvent ring). As before, we obtain conjectural density functions $\rho_{ref}^{geo}$ and $\rho_{ref}^{arith}$. Using our work in \S \ref{s:quartic}, we compute $\rho^{geo}$ and give the formula below.

\tpoint{Proposition} {\em
\label{p:quartic-rho}
In the notation of Section~\ref{s:quartic}, the function $\rho_{ref}^{geo}$ is supported on the union of $\qpolytope$ and 
\[
    \qpolytope' \cap \{(\ove_1,\ove_2,\ove_3;\ovf_1,\ovf_2) \in \RR^5 : \ovf_1 = 2\ove_1\}.
\]
On this polytope, we have
\[
    \rho_{ref}^{geo}(\ove_1,\ove_2,\ove_3;\ovf_1,\ovf_2) = 1 - \frac{1}{2}\sum_{ijk}\max\{0,\ovf_k - \ove_i - \ove_j\}.
\]
Moreover, the projection $\rho^{geo}$ is supported on the polytope $\mathcal{P}_4'$ and is given by the formula
\[
    \rho^{geo}(\ove_1,\ove_2,\ove_3) = \sup \{\rho_{ref}^{geo}(\ove_1,\ove_2,\ove_3;\ovf_1,\ovf_2) : \ovf_1 + \ovf_2 = 1\}.
\]
}

\bpf
Theorem~\ref{t:tet} classifies which $(e_1,e_2,e_3;f_1,f_2)$ arise from smooth quartic covers. When $\mathcal{H}_{\vec{e},\vec{f}}$ is nonempty, its dimension is given by
\[
    h^0(\proj^1, \Sym^2 \cE \otimes \cF^{\vee}) - h^0(\proj^1, \End(\oh(\vec{e}))) - h^0(\proj^1, \End(\oh(\vec{f}))) + 1.
\]
After expanding the expression above into a piecewise linear function in $\vec{e}$ and $\vec{f}$ and then taking appropriate limits, one obtains the claimed value for $\rho^{geo}_{ref}$. An analogous argument gives the claimed value for $\rho^{\geo}$.  We leave the details to the reader.
\epf

We may now make a conjecture  analogous to Conjecture~\ref{conj:rho}.

\tpoint{Conjecture}
\label{conj:rhoref}
{\em The quantities $\rho_{ref}^{geo}(x)$ and $\rho^{arith}_{ref}(x)$ exist and are equal. They satisfy 
\[
    \rho_{ref}^{geo} = \rho_{ref}^{arith}
\]
and are supported on a finite union of rational polytopes, piecewise linear, and continuous on their support.
}

\epoint{Final remarks: refined density functions in all degrees}
For $d = 5$, it is natural to make a conjecture analogous to Conjecture~\ref{conj:rhoref}. Canonically associated to a smooth degree $5$ cover is a degree $6$ cover called its sextic resolvent, and one should replace the data of $(f_1,f_2)$ in the quartic case with the data of the scrollar invariants of the sextic resolvent. Indeed from \cite[Theorem~1.5]{fv} one may compute an explicit description of $\rho_{ref}^{geo}$ and $\rho^{geo}$ for quintic covers. In this case, the density function $\rho_{ref}^{geo}$ is supported on a union of two polytopes, is piecewise linear, and is continuous on its support. 

For general degree $d$, the work of Casnati and Ekedahl \cite[Thm.~1.3]{casnatiekedahl}  (see also \cite{casnati1}) shows that a smooth degree $d$ cover $\pi : C \rightarrow \proj^1$ canonically embeds in the projectivization of the Tschirnhausen bundle. The minimal resolution of $C$ inside the $\PP^{d-2}$-bundle is of the form
\[
    0 \rightarrow \mathcal{N}_{d-2}(-d) \rightarrow \mathcal{N}_{d-3}(-d+2) \rightarrow \dots \rightarrow \mathcal{N}_1(-2) \rightarrow \oh_{\PP^1} \rightarrow \oh_C \rightarrow 0.
\]
The $\mathcal{N}_1,\dots,\mathcal{N}_{d-2}$ are vector bundles on $\PP^1$ whose ranks are determined by $d$, and whose degrees are of the form $c_i(d + g - 1)$ where $c_i$ is an integer depending only on $i$ and $d$.

In the case $d = 3$, the Tschirnhausen bundle $\cE_{\pi}$ determines every bundle in the minimal resolution. When $d = 4$ and $5$, the Tschirnhausen bundle of the associated cubic resolvent cover (resp. sextic resolvent cover), along with the Tschirnhausen bundle of the curve itself, determines every bundle in the minimal resolution. In general, one may consider the sublocus of $\cH_{d,g}$ where $\mathcal{N}_1,\dots,\mathcal{N}_{d-2}$ is fixed, and define a (conjectural) refined density function $\rho_{ref}^{geo}$. By using the tuning modules defined in \cite{woodyasuda}, which give lattices analogous to the vector bundles $\mathcal{N}_i$, one may define a (conjectural) arithmetic refined density function $\rho_{ref}^{arith}$. It is interesting to ask if the natural conjecture, analogous to Conjecture~\ref{conj:rhoref}, is true in this case.
}

\bpoint{Furthest extension: density functions valued in the Grothendieck ring}
There is a natural analogous question to consider in the Grothendieck ring of stacks, which we only state when $3 \leq d \leq 5$. Assume henceforth that $3 \leq d \leq 5$. 

Following \cite{ekedahl},
let $\grstk$ be the space obtained by taking the Grothendieck ring of varieties, inverting $\mathbb{L}$, and then completing with respect to the dimensional filtration.  Then $\grstk$ comes with a ''dimension function'' $\dim: \grstk \rightarrow \Z$, and $-\dim$ is a valuation on $\grstk$. 

\epoint{Motivation} { 
The goal of this section will be to define (conjecturally) a function $\rho^{gr}_{ref}$ taking values in $\grstk \otimes_{\Z} \R$. Moreover, if $\rho^{gr}_{ref}$ is defined, then  $\rho^{gr}_{ref}$ is $\rho^{geo}_{ref}(x)$  (a  real-valued function) plus negative-dimensional ``correction terms'':  $$\rho^{gr}_{ref}(x) = \rho^{geo}_{ref}(x) + Z(x)$$ for some function $Z(x)$ with $\dim(Z(x)) < 0$ for all $x$. (We consider $\rho^{geo}_{ref}(x)$ as taking values in $\grstk \otimes_\Z \R$ via the natural map $\RR \rightarrow \grstk \otimes_{\Z} \R$.) }


\epoint{Setup} {
Note first that $\pi^{geo}_{ref}$ is defined using \emph{quotients} of dimensions. To emulate this in $\grstk$, we will first partially define a logarithm with base $\mathbb{L}$ taking values in $\grstk \otimes_\Z \Q$, and take quotients of the form
\[
    \frac{\log_{\mathbb{L}}X}{\log_{\mathbb{L}}Y}
\]
where $X,Y \in \grstk$. If $\dim Y \neq 0$, then it will follow from our definition of logarithm (Definition~\ref{def:log}) that
\[
   \frac{\log_{\mathbb{L}}X}{\log_{\mathbb{L}}Y} = \frac{\dim X}{\dim Y} + Z
\]
where $\dim(Z) < 0$. Therefore, the quotient of the logarithms is the quotient of dimensions plus lower-dimensional terms, as desired. Conjecture~\ref{conj:rhoref} is about the convergence of $\dim X/\dim Y$, in an appropriate limit; in this section we extend to conjecture the convergence of the $Z$-term. 
}

\epoint{Definition} {\em
\label{def:log}
Let $X \in \grstk$. Suppose $X = \mathbb{L}^{\dim X} \times (1 + Y)$ where $\dim Y < 0$. Let $\log_{\mathbb{L}}(X) := \dim X + \sum_{i=1}^{\infty} (-1)^{i+1} \frac{Y^i}{i} \in \grstk \otimes_{\Z} \Q$.
}

As stated previously, to a degree $4$ cover there is a naturally associated pair $\vec{f} = (f_1,f_2)$, and to a degree $5$ cover there is a naturally associated quintuple $\vec{f} = (f_1,\dots,f_5)$. For the sake of convention, we associate to a degree $3$ cover the \emph{empty} tuple $\vec{f} = ()$.

\tpoint{Lemma} { \em
If $\mathcal{H}_{\vec{e},\vec{f}}$ is nonempty, then $[\mathcal{H}_{\vec{e},\vec{f}}] = \mathbb{L}^{\dim \mathcal{H}_{\vec{e},\vec{f}}} \times (1 + Y)$ for some $Y \in \grstk$ with $\dim(Y) < 0$. Similarly, $[\mathcal{H}_{d,g}] = \mathbb{L}^{\dim \mathcal{H}_{d,g}} \times (1 + Y)$ for some $Y \in \grstk$ with $\dim(Y) < 0$. Therefore $\log_{\mathbb{L}}[\mathcal{H}_{\vec{e},\vec{f}}]$ and $\log_{\mathbb{L}} [\mathcal{H}_{d,g}]$ are well-defined. Moreover $\log_{\mathbb{L}} [\mathcal{H}_{d,g}]$ is invertible in $\grstk \otimes_{\Z} \Q$. 
}

\bpf
We begin by proving that $[\mathcal{H}_{\vec{e},\vec{f}}]$ and $[\mathcal{H}_{d,g}]$ have the specified form. Each $\mathcal{H}_{d,g}$ has some $\mathcal{H}_{\vec{e},\vec{f}}$ as an open dense locus; therefore it suffices to prove the claim for $\mathcal{H}_{\vec{e},\vec{f}}$. By \cite{casnati, casnati1}, if $\mathcal{H}_{\vec{e},\vec{f}}$ is nonempty then
\[
    [\mathcal{H}_{\vec{e},\vec{f}}]\times [\GL(\oh(\vec{e}))] = \mathbb{L}^k \times (1 + X)
\]
\[
    [\mathcal{H}_{\vec{e},\vec{f}}]\times [\GL(\oh(\vec{e}))] \times [\GL(\oh(\vec{f}))]  = \mathbb{L}^k \times (1 + X)
\]
in the cases $d = 3$ and $d = 4,5$ respectively, for some integer $k$ and some element $X \in \grstk$ with $\dim(X) < 0$. A short computation shows that for any vector bundle $\oh(\vec{e})$ on $\proj^1$, we have $[\GL(\oh(\vec{e}))] = \LL^{\dim \GL(\oh(\vec{e}))}(1 + Z)$ for some $Z \in \grstk$ such that $\dim Z < 0$. The power series expansion
\begin{equation}
\label{e:pwr-series}
     \frac{1}{1 + z} = \sum_{n = 0}^{\infty}(-1)^nz^n.
\end{equation}
implies that $[\GL(\oh(\vec{e}))]^{-1} = \mathbb{L}^{-\dim \GL(\oh(\vec{e}))}\left (\sum_{n = 0}^{\infty}(-1)^nZ^n \right)$. Substituting, we obtain
\[
    [\mathcal{H}_{\vec{e},\vec{f}}] = \mathbb{L}^k \times (1 + X) \times [\GL(\oh(\vec{e}))]^{-1} = \mathbb{L}^k \times (1 + X) \times (1 - Z  + Z^2 \dots)
\]
\[
    [\mathcal{H}_{\vec{e},\vec{f}}] = \mathbb{L}^k \times (1 + X) \times [\GL(\oh(\vec{e}))]^{-1} \times [\GL(\oh(\vec{f}))]^{-1} = \mathbb{L}^k \times (1 + X) \times (1 - Z_e + \dots) \times (1 - Z_f + \dots),
\]
respectively in the cases $d = 3$ and $d = 4,5$ as required.

We now prove that $\log_{\mathbb{L}} [\mathcal{H}_{d,g}]$ is invertible in $\grstk \otimes_\Z \Q$. We have shown that
\[
    [\mathcal{H}_{d,g}] = \mathbb{L}^{\dim \mathcal{H}_{d,g}} \times (1 + Y)
\]
for some $Y \in \grstk$ with $\dim(Y) < 0$. By definition
\[
    \log_{\mathbb{L}} [\mathcal{H}_{d,g}] = \dim \mathcal{H}_{d,g} + \sum_{i = 1}^{\infty}(-1)^{i+1}\frac{Y^i}{i} = \dim \mathcal{H}_{d,g} \times (1 + W)
\]
where $W =  (\dim \mathcal{H}_{d,g})^{-1}\sum_{i = 1}^{\infty}(-1)^{i+1}\frac{Y^i}{i}$ so $\dim(W) < 0$. Therefore, $\log_{\mathbb{L}} [\mathcal{H}_{d,g}]^{-1}$ is given by
\[
    (\dim \mathcal{H}_{d,g})^{-1}(1 + W)^{-1} = (\dim \mathcal{H}_{d,g})^{-1}\sum_{n = 0}^{\infty}(-1)^nW^n \in \grstk \otimes_\Z \Q.
\]
\epf

When $3 \leq d \leq 5$ and $B$ is an open ball in $\RR^2$, $\RR^5$, or $\RR^9$ respectively, a potential alternate definition of $\pi$ is:
\[
    \pi_{ref, avg}^{geo}(B,g) = \frac{\avg \left \{ \dim \mathcal{H}_{\vec{e},\vec{f}} : (\vec{e},\vec{f})/(d+g-1) \in B, \; \mathcal{H}_{\vec{e},\vec{f}} \neq \emptyset \right\}}{\dim \mathcal{H}_{d,g}}.
\]
(By convention, we take $-1$ be the average of the empty set.) Although $\pi_{ref, avg}^{geo}(B,g)$ may not be equal to $\pi^{geo}_{ref}$, a short computation shows that the limit $\pi_{ref, avg}^{geo}(B)$ still exists, and moreover
\[
    \rho^{geo}_{ref} = \lim_{i \rightarrow \infty} \pi^{geo}_{ref}(B_i(x)) = \lim_{i \rightarrow \infty} \pi^{geo}_{ref,avg}(B_i(x)),
\]
where $\{B_i\}$ is any set of open balls all containing $x$ and with radius converging to $0$. In other words the resulting density functions are equal.     We make this equivalent alternate definition because it suggests a corresponding definition in the Grothendieck ring.

Thus  motivated, we define 
\[
    \pi^{gr}_{ref, avg}(B,g) = \frac{\avg\left\{{\log_{\mathbb{L}} [\mathcal{H}_{\vec{e},\vec{f}}]: (\vec{e},\vec{f})/(d+g-1) \in B, \; \mathcal{H}_{\vec{e},\vec{f}} \neq \emptyset}\right\}}{\log_{\mathbb{L}} [\mathcal{H}_{d,g}]}
\]
taking values in $\grstk \otimes_\Z \Q$. We may then proceed to make an analogous conjecture:

\tpoint{Conjecture (?)}
\label{conj:measuregr}
{\em Fix $3 \leq d \leq 5$. For all $B$,  $\lim_{g \rightarrow \infty} \pi^{gr}_{ref, avg}(B,g)$ converges in $\grstk \otimes_{\Z} \R$.}
\exercisedone

There is little direct evidence for this conjecture (perhaps it should be a question), but because  any counterexample to this conjecture would be very interesting, we feel it is worth making. As before, let $\pi^{gr}_{ref, avg}(B)$ denote the corresponding limit, and for a point $x$ in $\RR^2$, $\RR^5$, or $\RR^9$ respectively define
\[
    \rho^{gr}_{ref}(x) = \lim_{i \rightarrow \infty}\pi^{gr}_{ref,avg}(B_i(x))
\]
where $\{B_i(x)\}$ is any set of open balls all containing $x$ and with radius converging to $0$.

\tpoint{Conjecture}
\label{conj:rhogr} {\em
When $3 \leq d \leq 5$,  the limit $\rho^{gr}_{ref}(x)$ exists.
}




\end{document}

%% file: RVSVpreamble.tex
\usepackage{palatino, euler, epic,eepic, amssymb, xypic, floatflt, microtype, mathtools}
\usepackage[linktocpage,hidelinks]{hyperref}
\usepackage{verbatim,color}
\hypersetup{
    colorlinks,
    linkcolor={blue!50!black},
    citecolor={green!50!black},
    urlcolor={red!80!black}
}

\usepackage{tikz-cd}
\usepackage{url}
\usepackage{dsfont}

\input xy
\xyoption{all}

\newcommand{\bpf}{\noindent {\em Proof.  }}
\newcommand{\epf}{\qed \vspace{+10pt}}

\newcommand{\point}{\vspace{3mm}\par \noindent \refstepcounter{subsection}{\thesubsection.} }
\newcommand{\tpoint}[1]{\vspace{3mm}\par \noindent \refstepcounter{subsection}{\thesubsection.} 
  {\bf #1. ---} }
\newcommand{\epoint}[1]{\vspace{3mm}\par \noindent \refstepcounter{subsection}{\thesubsection.} 
  {\em #1.} }
\newcommand{\bpoint}[1]{\vspace{3mm}\par \noindent \refstepcounter{subsection}{\thesubsection.} 
  {\bf #1.} }

\newcommand{\exercisedone}{ \vspace{2mm}}

\newcommand{\proj}{\mathbb P}

\newcommand{\C}{\mathbb{C}}
\newcommand{\F}{\mathbb{F}}

\newcommand{\Z}{\mathbb{Z}}

\newcommand{\Q}{\mathbb{Q}}

\newcommand{\R}{\mathbb{R}}
\def\CC{\mathbb{C}}

\def\FF{\mathbb{F}}

\def\LL{\mathbb{L}}

\def\PP{\mathbb{P}}
\def\QQ{\mathbb{Q}}
\def\RR{\mathbb{R}}

\def\ZZ{\mathbb{Z}}

\newcommand{\cA}{{\mathscr{A}}}

\newcommand{\cE}{{\mathscr{E}}}
\newcommand{\cF}{{\mathscr{F}}}

\newcommand{\oh}{{\mathscr{O}}}

\newcommand{\Ga}{\Gamma}

\newcommand{\geo}{\rm{geo}}

\newcommand{\coloneq}{:=}

\newcommand{\propernormal}{%
  \mathrel{\ooalign{$\lneq$\cr\raise.22ex\hbox{$\lhd$}\cr}}}

\newcommand\Stab{\operatorname{Stab}}

\newcommand{\Sym}{\operatorname{Sym}}

\newcommand{\Hom}{\operatorname{Hom}}

\newcommand{\Pic}{\operatorname{Pic}}
\newcommand{\Spec}{\operatorname{Spec}}

\newcommand{\Stck}{\operatorname{Stck}}

\newcommand{\lremind}[1]{{\bf[label:  #1]}}
\newcommand{\notation}[1]{}
\renewcommand{\lremind}[1]{{}}

\newcommand{\cut}[1]{}

\DeclareMathOperator{\GL}{GL}

\DeclareMathOperator{\End}{End}

\DeclareMathOperator{\avg}{avg}

\DeclareMathOperator{\Disc}{Disc}

\DeclareMathOperator{\Proj}{Proj}

\DeclareMathOperator{\grstk}{\widehat{K(\Stck)}}

\newcommand{\Hur}{\mathcal{H}ur}    
\newcommand{\Bun}{\mathcal{B}un}    
\newcommand{\AAA}{\mathbb{A}}

\newcommand{\polytope}{\mathcal{P}}  
\newcommand{\qpolytope}{\mathcal{Q}}
\newcommand{\ove}{\overline{e}}
\newcommand{\ovf}{\overline{f}}

\renewcommand{\cH}{\mathcal{H}}

\newcommand{\Tsch}{\operatorname{Tsch}}